\documentclass{amsart}
\usepackage{amssymb}

\newtheorem{theorem}{Theorem}[section]
\newtheorem{corollary}[theorem]{Corollary}
\newtheorem{conjecture}[theorem]{Conjecture}
\newtheorem{lemma}[theorem]{Lemma}
\newtheorem{proposition}[theorem]{Proposition}
\newtheorem{definition}[theorem]{Definition}
\newtheorem{axiom}{Axiom}[section]
\newtheorem{remark}{Remark}[section]
\newtheorem{example}{Example}[section]
\newtheorem{exercise}{Exercise}[section]

\numberwithin{equation}{section}

\def\forces{\mathbin{\parallel\mkern-9mu-}}
\makeatletter
\newcount\@hour\newcount\@minute\chardef\@x10\chardef\@xv60
\def\tcitime{
\def\@time{%
  \@minute\time\@hour\@minute\divide\@hour\@xv
  \ifnum\@hour<\@x 0\fi\the\@hour:%
  \multiply\@hour\@xv\advance\@minute-\@hour
  \ifnum\@minute<\@x 0\fi\the\@minute
  }}%

\@ifundefined{hyperref}{}{}

\@ifundefined{qExtProgCall}{\def\qExtProgCall#1#2#3#4#5#6{\relax}}{}
\def\QCTOpt[#1]#2{%
  \def\QCTOptB{#1}
  \def\QCTOptA{#2}
}
\def\QCTNOpt#1{%
  \def\QCTOptA{#1}
  \let\QCTOptB\empty
}
\def\Qct{%
  \@ifnextchar[{%
    \QCTOpt}{\QCTNOpt}
}
\def\QCBOpt[#1]#2{%
  \def\QCBOptB{#1}
  \def\QCBOptA{#2}
}
\def\QCBNOpt#1{%
  \def\QCBOptA{#1}
  \let\QCBOptB\empty
}
\def\Qcb{%
  \@ifnextchar[{%
    \QCBOpt}{\QCBNOpt}
}
\def\PrepCapArgs{%
  \ifx\QCBOptA\empty
    \ifx\QCTOptA\empty
      {}%
    \else
      \ifx\QCTOptB\empty
        {\QCTOptA}%
      \else
        [\QCTOptB]{\QCTOptA}%
      \fi
    \fi
  \else
    \ifx\QCBOptA\empty
      {}%
    \else
      \ifx\QCBOptB\empty
        {\QCBOptA}%
      \else
        [\QCBOptB]{\QCBOptA}%
      \fi
    \fi
  \fi
}
\newcount\GRAPHICSTYPE
\GRAPHICSTYPE=\z@
\def\GRAPHICSPS#1{%
 \ifcase\GRAPHICSTYPE
   \special{ps: #1}%
 \or
   \special{language "PS", include "#1"}%
 \fi
}%
\def\graffile#1#2#3#4{%
    \leavevmode
    \raise -#4 \BOXTHEFRAME{%
        \hbox to #2{\raise #3\hbox to #2{\null #1\hfil}}}%
}%
\def\draftbox#1#2#3#4{%
 \leavevmode\raise -#4 \hbox{%
  \frame{\rlap{\protect\tiny #1}\hbox to #2%
   {\vrule height#3 width\z@ depth\z@\hfil}%
  }%
 }%
}%
\newcount\draft
\draft=\z@

\newif\ifwasdraft
\wasdraftfalse

\def\GRAPHIC#1#2#3#4#5{%
 \ifnum\draft=\@ne\draftbox{#2}{#3}{#4}{#5}%
  \else\graffile{#1}{#3}{#4}{#5}%
  \fi
 }%
\def\addtoLaTeXparams#1{%
    \edef\LaTeXparams{\LaTeXparams #1}}%
\newif\ifBoxFrame \BoxFramefalse
\newif\ifOverFrame \OverFramefalse
\newif\ifUnderFrame \UnderFramefalse

\def\BOXTHEFRAME#1{%
   \hbox{%
      \ifBoxFrame
         \frame{#1}%
      \else
         {#1}%
      \fi
   }%
}

\def\doFRAMEparams#1{\BoxFramefalse\OverFramefalse\UnderFramefalse\readFRAMEparams#1\end}%
\def\readFRAMEparams#1{%
 \ifx#1\end%
  \let\next=\relax
  \else
  \ifx#1i\dispkind=\z@\fi
  \ifx#1d\dispkind=\@ne\fi
  \ifx#1f\dispkind=\tw@\fi
  \ifx#1t\addtoLaTeXparams{t}\fi
  \ifx#1b\addtoLaTeXparams{b}\fi
  \ifx#1p\addtoLaTeXparams{p}\fi
  \ifx#1h\addtoLaTeXparams{h}\fi
  \ifx#1X\BoxFrametrue\fi
  \ifx#1O\OverFrametrue\fi
  \ifx#1U\UnderFrametrue\fi
  \ifx#1w
    \ifnum\draft=1\wasdrafttrue\else\wasdraftfalse\fi
    \draft=\@ne
  \fi
  \let\next=\readFRAMEparams
  \fi
 \next
 }%
\def\IFRAME#1#2#3#4#5#6{%
      \bgroup
      \let\QCTOptA\empty
      \let\QCTOptB\empty
      \let\QCBOptA\empty
      \let\QCBOptB\empty
      #6%
      \parindent=0pt%
      \leftskip=0pt
      \rightskip=0pt
      \setbox0 = \hbox{\QCBOptA}%
      \@tempdima = #1\relax
      \ifOverFrame
          \typeout{This is not implemented yet}%
          \show\HELP
      \else
         \ifdim\wd0>\@tempdima
            \advance\@tempdima by \@tempdima
            \ifdim\wd0 >\@tempdima
               \textwidth=\@tempdima
               \setbox1 =\vbox{%
                  \noindent\hbox to \@tempdima{\hfill\GRAPHIC{#5}{#4}{#1}{#2}{#3}\hfill}\\%
                  \noindent\hbox to \@tempdima{\parbox[b]{\@tempdima}{\QCBOptA}}%
               }%
               \wd1=\@tempdima
            \else
               \textwidth=\wd0
               \setbox1 =\vbox{%
                 \noindent\hbox to \wd0{\hfill\GRAPHIC{#5}{#4}{#1}{#2}{#3}\hfill}\\%
                 \noindent\hbox{\QCBOptA}%
               }%
               \wd1=\wd0
            \fi
         \else
            \ifdim\wd0>0pt
              \hsize=\@tempdima
              \setbox1 =\vbox{%
                \unskip\GRAPHIC{#5}{#4}{#1}{#2}{0pt}%
                \break
                \unskip\hbox to \@tempdima{\hfill \QCBOptA\hfill}%
              }%
              \wd1=\@tempdima
           \else
              \hsize=\@tempdima
              \setbox1 =\vbox{%
                \unskip\GRAPHIC{#5}{#4}{#1}{#2}{0pt}%
              }%
              \wd1=\@tempdima
           \fi
         \fi
         \@tempdimb=\ht1
         \advance\@tempdimb by \dp1
         \advance\@tempdimb by -#2%
         \advance\@tempdimb by #3%
         \leavevmode
         \raise -\@tempdimb \hbox{\box1}%
      \fi
      \egroup%
}%
\def\DFRAME#1#2#3#4#5{%
 \begin{center}
     \let\QCTOptA\empty
     \let\QCTOptB\empty
     \let\QCBOptA\empty
     \let\QCBOptB\empty
     \ifOverFrame 
        #5\QCTOptA\par
     \fi
     \GRAPHIC{#4}{#3}{#1}{#2}{\z@}
     \ifUnderFrame 
        \nobreak\par #5\QCBOptA
     \fi
 \end{center}%
 }%
\def\FFRAME#1#2#3#4#5#6#7{%
 \begin{figure}[#1]%
  \let\QCTOptA\empty
  \let\QCTOptB\empty
  \let\QCBOptA\empty
  \let\QCBOptB\empty
  \ifOverFrame
    #4
    \ifx\QCTOptA\empty
    \else
      \ifx\QCTOptB\empty
        \caption{\QCTOptA}%
      \else
        \caption[\QCTOptB]{\QCTOptA}%
      \fi
    \fi
    \ifUnderFrame\else
      \label{#5}%
    \fi
  \else
    \UnderFrametrue%
  \fi
  \begin{center}\GRAPHIC{#7}{#6}{#2}{#3}{\z@}\end{center}%
  \ifUnderFrame
    #4
    \ifx\QCBOptA\empty
      \caption{}%
    \else
      \ifx\QCBOptB\empty
        \caption{\QCBOptA}%
      \else
        \caption[\QCBOptB]{\QCBOptA}%
      \fi
    \fi
    \label{#5}%
  \fi
  \end{figure}%
 }%
\newcount\dispkind%

\def\makeactives{
  \catcode`\"=\active
  \catcode`\;=\active
  \catcode`\:=\active
  \catcode`\'=\active
  \catcode`\~=\active
}
\bgroup
   \makeactives
   \gdef\activesoff{%
      \def"{\string"}
      \def;{\string;}
      \def:{\string:}
      \def'{\string'}
      \def~{\string~}
    }
\egroup

\def\FRAME#1#2#3#4#5#6#7#8{%
 \bgroup
 \@ifundefined{bbl@deactivate}{}{\activesoff}
 \ifnum\draft=\@ne
   \wasdrafttrue
 \else
   \wasdraftfalse%
 \fi
 \def\LaTeXparams{}%
 \dispkind=\z@
 \def\LaTeXparams{}%
 \doFRAMEparams{#1}%
 \ifnum\dispkind=\z@\IFRAME{#2}{#3}{#4}{#7}{#8}{#5}\else
  \ifnum\dispkind=\@ne\DFRAME{#2}{#3}{#7}{#8}{#5}\else
   \ifnum\dispkind=\tw@
    \edef\@tempa{\noexpand\FFRAME{\LaTeXparams}}%
    \@tempa{#2}{#3}{#5}{#6}{#7}{#8}%
    \fi
   \fi
  \fi
  \ifwasdraft\draft=1\else\draft=0\fi{}%
  \egroup
 }%
\def\TEXUX#1{"texux"}

\def\limfunc#1{\mathop{\rm #1}}%

\long\def\QQQ#1#2{%
     \long\expandafter\def\csname#1\endcsname{#2}}%
\@ifundefined{QTP}{\def\QTP#1{}}{}
\@ifundefined{QEXCLUDE}{\def\QEXCLUDE#1{}}{}
\@ifundefined{Qlb}{}{}
\@ifundefined{Qlt}{}{}
\long\def\QQA#1#2{}%
\def\QTR#1#2{{\csname#1\endcsname #2}}
\def\EXPAND#1[#2]#3{}%
\def\NOEXPAND#1[#2]#3{}%
\def\LaTeXparent#1{}%
\def\ChildStyles#1{}%
\def\ChildDefaults#1{}%
\def\QTagDef#1#2#3{}%
\@ifundefined{StyleEditBeginDoc}{}{}
\def\QQfnmark#1{\footnotemark}

\@ifundefined{INDEX}{\def\INDEX#1#2{}{}}{}%
\@ifundefined{SUBINDEX}{\def\SUBINDEX#1#2#3{}{}{}}{}%
\@ifundefined{initial}%
   {\def\initial#1{\bigbreak{\raggedright\large\bf #1}\kern 2\p@\penalty3000}}%
   {}%
\@ifundefined{entry}{}{}%
\@ifundefined{primary}{}{}%
\@ifundefined{secondary}{}{}%
\@ifundefined{ZZZ}{}{\makeatletter\input gnuindex.sty\makeatother\makeindex\makeatletter}%
\@ifundefined{abstract}{%
 \def\abstract{%
  \if@twocolumn
   \section*{Abstract (Not appropriate in this style!)}%
   \else \small 
   \begin{center}{\bf Abstract\vspace{-.5em}\vspace{\z@}}\end{center}%
   \quotation 
   \fi
  }%
 }{%
 }%
\@ifundefined{endabstract}{\def\endabstract
  {\if@twocolumn\else\endquotation\fi}}{}%
\@ifundefined{maketitle}{\def\maketitle#1{}}{}%
\@ifundefined{affiliation}{\def\affiliation#1{}}{}%
\@ifundefined{proof}{}{}%
\@ifundefined{endproof}{}{}%
\@ifundefined{newfield}{\def\newfield#1#2{}}{}%
\@ifundefined{chapter}{\def\chapter#1{\par(Chapter head:)#1\par }%
 \newcount\c@chapter}{}%
\@ifundefined{part}{\def\part#1{\par(Part head:)#1\par }}{}%
\@ifundefined{section}{\def\section#1{\par(Section head:)#1\par }}{}%
\@ifundefined{subsection}{\def\subsection#1%
 {\par(Subsection head:)#1\par }}{}%
\@ifundefined{subsubsection}{\def\subsubsection#1%
 {\par(Subsubsection head:)#1\par }}{}%
\@ifundefined{paragraph}{\def\paragraph#1%
 {\par(Subsubsubsection head:)#1\par }}{}%
\@ifundefined{subparagraph}{\def\subparagraph#1%
 {\par(Subsubsubsubsection head:)#1\par }}{}%
\@ifundefined{therefore}{}{}%
\@ifundefined{backepsilon}{}{}%
\@ifundefined{yen}{}{}%
\@ifundefined{registered}{%
   \def\registered{\relax\ifmmode{}\r@gistered
                    \else$\m@th\r@gistered$\fi}%
 \def\r@gistered{^{\ooalign
  {\hfil\raise.07ex\hbox{$\scriptstyle\rm\text{R}$}\hfil\crcr
  \mathhexbox20D}}}}{}%
\@ifundefined{Eth}{}{}%
\@ifundefined{eth}{}{}%
\@ifundefined{Thorn}{}{}%
\@ifundefined{thorn}{}{}%
\@ifundefined{degree}{}{}%
\newdimen\theight
\def\Column{%
 \vadjust{\setbox\z@=\hbox{\scriptsize\quad\quad tcol}%
  \theight=\ht\z@\advance\theight by \dp\z@\advance\theight by \lineskip
  \kern -\theight \vbox to \theight{%
   \rightline{\rlap{\box\z@}}%
   \vss
   }%
  }%
 }%
\def\qed{%
 \ifhmode\unskip\nobreak\fi\ifmmode\ifinner\else\hskip5\p@\fi\fi
 \hbox{\hskip5\p@\vrule width4\p@ height6\p@ depth1.5\p@\hskip\p@}%
 }%
\def\miss{\hbox{\vrule height2\p@ width 2\p@ depth\z@}}%
\def\tcol#1{{\baselineskip=6\p@ \vcenter{#1}} \Column}  %
\def\newfmtname{LaTeX2e}
\def\chkcompat{%
   \if@compatibility
   \else
     \usepackage{latexsym}
   \fi
}

\ifx\fmtname\newfmtname
  \DeclareOldFontCommand{\rm}{\normalfont\rmfamily}{\mathrm}
  \DeclareOldFontCommand{\sf}{\normalfont\sffamily}{\mathsf}
  \DeclareOldFontCommand{\tt}{\normalfont\ttfamily}{\mathtt}
  \DeclareOldFontCommand{\bf}{\normalfont\bfseries}{\mathbf}
  \DeclareOldFontCommand{\it}{\normalfont\itshape}{\mathit}
  \DeclareOldFontCommand{\sl}{\normalfont\slshape}{\@nomath\sl}
  \DeclareOldFontCommand{\sc}{\normalfont\scshape}{\@nomath\sc}
  \chkcompat
\fi

\def\alpha{{\Greekmath 010B}}%
\def\beta{{\Greekmath 010C}}%
\def\gamma{{\Greekmath 010D}}%
\def\delta{{\Greekmath 010E}}%
\def\epsilon{{\Greekmath 010F}}%
\def\zeta{{\Greekmath 0110}}%
\def\eta{{\Greekmath 0111}}%
\def\theta{{\Greekmath 0112}}%
\def\iota{{\Greekmath 0113}}%
\def\kappa{{\Greekmath 0114}}%
\def\lambda{{\Greekmath 0115}}%
\def\mu{{\Greekmath 0116}}%
\def\nu{{\Greekmath 0117}}%
\def\xi{{\Greekmath 0118}}%
\def\pi{{\Greekmath 0119}}%
\def\rho{{\Greekmath 011A}}%
\def\sigma{{\Greekmath 011B}}%
\def\tau{{\Greekmath 011C}}%
\def\upsilon{{\Greekmath 011D}}%
\def\phi{{\Greekmath 011E}}%
\def\chi{{\Greekmath 011F}}%
\def\psi{{\Greekmath 0120}}%
\def\omega{{\Greekmath 0121}}%
\def\varepsilon{{\Greekmath 0122}}%
\def\vartheta{{\Greekmath 0123}}%
\def\varpi{{\Greekmath 0124}}%
\def\varrho{{\Greekmath 0125}}%
\def\varsigma{{\Greekmath 0126}}%
\def\varphi{{\Greekmath 0127}}%

\def\nabla{{\Greekmath 0272}}
\def\FindBoldGroup{%
   {\setbox0=\hbox{$\mathbf{x\global\edef\theboldgroup{\the\mathgroup}}$}}%
}

\def\Greekmath#1#2#3#4{%
    \if@compatibility
        \ifnum\mathgroup=\symbold
           \mathchoice{\mbox{\boldmath$\displaystyle\mathchar"#1#2#3#4$}}%
                      {\mbox{\boldmath$\textstyle\mathchar"#1#2#3#4$}}%
                      {\mbox{\boldmath$\scriptstyle\mathchar"#1#2#3#4$}}%
                      {\mbox{\boldmath$\scriptscriptstyle\mathchar"#1#2#3#4$}}%
        \else
           \mathchar"#1#2#3#4%
        \fi 
    \else 
        \FindBoldGroup
        \ifnum\mathgroup=\theboldgroup 
           \mathchoice{\mbox{\boldmath$\displaystyle\mathchar"#1#2#3#4$}}%
                      {\mbox{\boldmath$\textstyle\mathchar"#1#2#3#4$}}%
                      {\mbox{\boldmath$\scriptstyle\mathchar"#1#2#3#4$}}%
                      {\mbox{\boldmath$\scriptscriptstyle\mathchar"#1#2#3#4$}}%
        \else
           \mathchar"#1#2#3#4%
        \fi                 
          \fi}

\newif\ifGreekBold  \GreekBoldfalse
\let\SAVEPBF=\pbf
\def\pbf{\GreekBoldtrue\SAVEPBF}%

\@ifundefined{theorem}{\newtheorem{theorem}{Theorem}}{}
\@ifundefined{lemma}{\newtheorem{lemma}[theorem]{Lemma}}{}
\@ifundefined{corollary}{\newtheorem{corollary}[theorem]{Corollary}}{}
\@ifundefined{conjecture}{}{}
\@ifundefined{proposition}{\newtheorem{proposition}[theorem]{Proposition}}{}
\@ifundefined{axiom}{}{}
\@ifundefined{remark}{}{}
\@ifundefined{example}{}{}
\@ifundefined{exercise}{}{}
\@ifundefined{definition}{\newtheorem{definition}{Definition}}{}

\@ifundefined{mathletters}{%
  \newcounter{equationnumber}  
  \def\mathletters{%
     \addtocounter{equation}{1}
     \edef\@currentlabel{\theequation}%
     \setcounter{equationnumber}{\c@equation}
     \setcounter{equation}{0}%
     \edef\theequation{\@currentlabel\noexpand\alph{equation}}%
  }
  
}{}

\@ifundefined{BibTeX}{%
    \def\BibTeX{{\rm B\kern-.05em{\sc i\kern-.025em b}\kern-.08em
                 T\kern-.1667em\lower.7ex\hbox{E}\kern-.125emX}}}{}%
\@ifundefined{AmS}%
    {\def\AmS{{\protect\usefont{OMS}{cmsy}{m}{n}%
                A\kern-.1667em\lower.5ex\hbox{M}\kern-.125emS}}}{}%
\@ifundefined{AmSTeX}{}{}%
\begin{document}
\title[A Non-reflexive Whitehead Group]{A Non-reflexive Whitehead Group}
\author{Paul C. Eklof}
\address[Eklof]{Math Dept, UCI\\
Irvine, CA 92697-3875}
\email{peklof@math.uci.edu}
\thanks{Partially supported by NSF Grants DMS-9501415 and DMS-9704477}
\author{Saharon Shelah}
\address[Shelah]{Institute of Mathematics, Hebrew University\\
Jerusalem 91904, Israel}
\email{shelah@math.huji.ac.il}
\thanks{Research of second author supported by German-Israeli Foundation for
Scientific Research \& Development Grant No. G-294.081.06/93. Pub. No. 621}
\date{}

\begin{abstract}
We prove that it is consistent that there is a non-reflexive Whitehead
group, in fact one whose dual group is free. We also prove that it is
consistent that there is a group $A$ such that $\limfunc{Ext}(A,{\mathbb Z})$
is torsion and $\limfunc{Hom}(A,{\mathbb Z})=0$. As an application we show
the consistency of the existence of new co-Moore spaces.
\end{abstract}

\maketitle

\section{Introduction}

This paper is motivated by a theorem and a question due to Martin Huber. He
proved \cite{Hu} in ZFC that if $A$ is $\aleph _{1}$-coseparable (that is, $%
\limfunc{Ext}(A,\mathbb{Z}^{(\omega )})=0$), then $A$ is reflexive (that is,
the natural map of $A$ to its double dual $A^{**}=\limfunc{Hom}(\limfunc{Hom}%
(A,\mathbb{Z}),\mathbb{Z})$ is an isomorphism). He asked whether it is provable in
ZFC that every Whitehead group $A$ (i.e., $\limfunc{Ext}(A,\mathbb{Z})=0$) is
reflexive. This is true in any model where every Whitehead group is free. It
is also true for Whitehead groups of cardinality $\aleph _{1}$ in a model of
MA + $\lnot $CH (because they are $\aleph _{1}$-coseparable: cf. 
\cite[Cor. XII.1.12]{EM}). Moreover, it is true in the original models of
GCH where there are non-free Whitehead groups (cf. \cite{Sh77}, \cite{Sh80}, 
\cite[Thm. XII.1.9]{EM}). It was left as an open question in \cite[p. 455]
{EM} whether every Whitehead group is reflexive. Here we give a strong
negative answer:

\begin{theorem}
\label{nonref}It is consistent with ZFC that there is a strongly
non-reflexive strongly $\aleph _{1}$-free Whitehead group $A$ of cardinality 
$\aleph _{1}$.
\end{theorem}

A group $A$ is strongly non-reflexive if $A$ is not isomorphic to $A^{**}$.
In fact, the example $A$ has the property that $A^{*}$ is free of rank $%
\aleph _{2}$ (i.e., isomorphic to $\mathbb{Z}^{(\aleph _{2})}$) so $A^{**}$ is
isomorphic to the product $\mathbb{Z}^{\aleph _{2}}$; it is therefore not 
isomorphic to $A$ since its cardinality is $2^{\aleph _{2}}>\aleph _{1}$.
(See Theorem \ref{free} and Corollary \ref{freecor} of section 1.)

If $\limfunc{Ext}(A,\mathbb{Z})=0$, then $A$ is separable (\cite[Thm 99.1]{F})
and hence $A^{*}$ is non-zero. However, using the same methods we can also
prove:

\begin{theorem}
\label{torsion}It is consistent with ZFC that there is a non-free strongly $%
\aleph _{1}$-free group $A$ of cardinality $\aleph _{1}$ such that $\limfunc{%
Ext}(A,\mathbb{Z})$ is torsion and $\limfunc{Hom}(A,\mathbb{Z})=0$.
\end{theorem}

It is not a theorem of ZFC that there is a non-free torsion-free group 
$A$  such
that $\limfunc{Ext}(A,\mathbb{Z})$ is torsion. Indeed, in any model where every
Whitehead group is free---a hypothesis which is consistent with CH or $\lnot $%
CH (cf. \cite{MS})---if  $A$ is not free, then $\limfunc{Ext}(A,\mathbb{Z}\mathbf{%
)}$ is not torsion (\cite{HHS}, \cite{EH}, \cite[ Thm. XII.2.4]{EM}).

Theorems \ref{free} and \ref{torsion} provide new examples of possible
co-Moore spaces (see section 6). In particular, we answer a question in 
\cite[p. 46]{GG} by showing that it is consistent that for any 
$n\geq 2$ there is a
co-Moore space of type $(F,n)$ where $F$ is a free group of rank $\aleph _{2}
$.

\smallskip\ 

The models for both theorems result from a finite support iteration of
c.c.c. posets and are models of ZFC + $\lnot $CH. (Other methods will be
needed to obtain consistency with CH.) We begin the iteration with a poset
which yields ``generic data'' from which the group $A$ is defined; we then
iterate the natural posets which insure that $\limfunc{Ext}(A,\mathbb{Z})=0$
(resp. $\limfunc{Ext}(A,\mathbb{Z})$ is torsion). The hard work is in proving
that $\limfunc{Hom}(A,\mathbb{Z})$ is as claimed. We define the forcing and the
group more precisely in the next section and then prove their properties in
the succeeding sections.

\section{The basic construction}

The group-theoretic construction is a generalization of that in 
\cite[XII.3.4]{EM}. Let $E$ be a stationary and co-stationary subset of $%
\omega _{1}$ consisting of limit ordinals, and for each $\delta \in E$, let $%
\eta _{\delta }$ be a ladder on $\delta $, that is, a strictly increasing
function $\eta _{\delta }:\omega \rightarrow \delta $ whose range approaches 
$\delta $. Let $F$ be the free abelian group with basis $\{x_{\nu }\colon %
\nu \in \omega _{1}\}\cup \{z_{\delta ,n}\colon \delta \in E,n\in \omega \}$%
. Let $g$ be a function from $E\times \omega $ to the integers $\geq 1$. Let 
$u$ be a function from $E\times \omega $ to the subgroup $\left\langle
x_{\nu }\colon \nu \in \omega _{1}\right\rangle $ generated by $\{x_{\nu }%
\colon \nu \in \omega _{1}\}$ such that $u(\delta ,n)$ belongs to $%
\left\langle x_{\nu }\colon \nu <\eta _{\delta }(n)\right\rangle $. Let $K$
be the subgroup of $F$ generated by $\{w_{\delta ,n}:\delta \in E,n\in
\omega \}$ where 
\begin{equation}
\text{ }w_{\delta ,n}=2^{g(\delta ,n)}z_{\delta ,n+1}-z_{\delta ,n}-x_{\eta
_{\delta }(n)}-u(\delta ,n).  \label{1.1}
\end{equation}
Let $A=F/K$. Then clearly $A$ is an abelian group of cardinality $\aleph _{1}
$. Notice that because the right-hand side of (\ref{1.1}) is 0 in $A$, we
have for each $\delta \in E$ and $n\in \omega $ the following relations in $A
$: 
\begin{equation}
2^{g(\delta ,n)}z_{\delta ,n+1}=z_{\delta ,n}+x_{\eta _{\delta
}(n)}+u(\delta ,n)  \label{1.15}
\end{equation}
and 
\begin{equation}
\text{ }2^{\sum_{j=0}^{n}g(\delta ,j)}z_{_{\delta ,n+1}}=z_{\delta
,0}+\sum_{k=0}^{n}\text{ }2^{\sum_{j=0}^{k-1}g(\delta ,j)}(x_{\eta _{\delta
}(k)}+u(\delta ,k))  \label{1.2}
\end{equation}
Here, and occasionally in what follows, we abuse notation and write, for
example, $z_{_{\delta ,n+1}}$ instead of $z_{_{\delta ,n+1}}+K$ for an
element of $A$. For each $\alpha <\omega _{1}$, let $A_{\alpha }$ be the
subgroup of $A$ generated by 
\begin{equation}
\{x_{\nu }:\nu <\alpha \}\cup \{z_{\delta ,n}:\delta \in E\cap \alpha \text{%
, }n\in \omega \}\text{.}  \label{1.4}
\end{equation}
Then, by (\ref{1.2}), for each $\delta \in E$, $z_{\delta ,0}+A_{\delta }$
is non-zero and divisible in $A_{\delta +1}/A_{\delta }$ by $2^{m}$ for all $%
m\in \omega $. Thus $A_{\delta +1}/A_{\delta }$ is not free and hence $A$ is
not free. (In fact $\Gamma (A)\supseteq \tilde{E}$.) Moreover, $A$ is
strongly $\aleph _{1}$-free; in fact, for every $\alpha <\omega _{1}$, using
Pontryagin's Criterion \cite[IV.2.3]{EM} we can show that $A/A_{\alpha +1}$
is $\aleph _{1}$-free for all $\alpha \in \omega _{1}\cup \{-1\}$.

We begin with a model $V$ of ZFC where GCH holds, choose $E\in V$, and
define the group $A$ in a generic extension $V^{Q_{0}}$ using generic
ladders $\eta _{\delta }$, and generic $u$ and $g$. Specifically:

\begin{definition}
\label{Q0}Let $Q_{0}$ be the set of all finite functions $q$ such that $%
\limfunc{dom}(q)$ is a finite subset of $E$ and for all $\gamma \in \limfunc{%
dom}(q)$, $q(\gamma )$ is a triple $(\eta _{\gamma }^{q},u_{\gamma
}^{q},g_{\gamma }^{q})$ where for some $r_{\gamma }^{q}\in \omega $:

\begin{itemize}
\item  $\eta _{\gamma }^{q}$ is a strictly increasing function$:r_{\gamma
}^{q}\rightarrow \gamma $ ;

\item  $u_{\gamma }^{q}:\{\gamma \}\times r_{\gamma }^{q}\rightarrow
\left\langle x_{\nu }\colon \nu \in \omega _{1}\right\rangle $ such that for
all $n<r_{\gamma }^{q}$, $u_{\gamma }^{q}(\gamma ,n)\in \left\langle x_{\nu
}:\nu <\eta _{\gamma }(n)\right\rangle $; and

\item  $g_{\gamma }^{q}:\{\gamma \}\times r_{\gamma }^{q}\rightarrow \{n\in
\omega :n\geq 1\}$.
\end{itemize}
\end{definition}

The partial ordering is defined by: $q_{1}\leq q_{2}$ if and only if $%
q_{1}\subseteq q_{2}$; note that we follow the convention that stronger
conditions are larger. Clearly $Q_{0}$ is c.c.c. and hence $E$ remains
stationary and co-stationary in a generic extension. We now do an iterated
forcing to make $A$ a Whitehead group. We begin by defining the basic
forcing that we will iterate.

\begin{definition}
Given a homomorphism $\psi :K\rightarrow \mathbb{Z}$, let $Q_{\psi }$ be the
poset of all finite functions $q$ into $\mathbb{Z}$ satisfying:

There are $\delta _{0}<\delta _{1}<...<\delta _{m}$ in $E$ and $\{r_{\ell
}:\ell \leq m\}\subseteq \omega $ such that $\limfunc{dom}(q)=$%
\[
\{z_{\delta _{\ell },n}:\ell \leq m,n\leq r_{\ell }\}\cup \{x_{\nu }:\nu \in
I_{q}\}
\]
where $I_{q}\subset \omega _{1}$ is finite and is such that for all $\ell
\leq m$

\begin{equation}
n<r_{\ell }\Rightarrow u(\delta _{\ell },n)\in \left\langle x_{\nu }:\nu \in
I_{q}\right\rangle \text{ and }\eta _{\delta _{\ell }}(n)\in
I_{q}\Leftrightarrow n<r_{\ell }  \label{5}
\end{equation}

\noindent and for all $\ell \leq m$ and $n<r_{\ell }$, $u(\delta _{\ell
},n)\in \left\langle x_{\nu }:\nu \in I_{q}\right\rangle $ and

\begin{equation}
\psi (w_{\delta _{\ell },n})=\text{ }2^{g(\delta ,n)}q(z_{\delta _{\ell
},n+1})-q(z_{\delta _{\ell },n})-q(x_{\eta _{\delta _{\ell
}}(n)})-q(u(\delta _{\ell },n)).  \label{6}
\end{equation}

\noindent (Compare with (\ref{1.1}). The definition of $q(u(\delta _{\ell
},n))$ is the obvious one, given that $q$ should extend to a homomorphism.)
Moreover, we require of $q$ that for all $\ell \neq j$ in $\{0,...,m\}$, 
\begin{equation}
\text{ }\eta _{\delta _{j}}(k)\neq \eta _{\delta _{\ell }}(i)\text{ for all }%
k\geq r_{j}\text{ and }i\in \omega .  \label{7}
\end{equation}
\end{definition}

We will denote $\{\delta _{0},...,\delta _{m}\}$ by $\limfunc{cont}(q)$ and $%
r_{\ell }$ by $\limfunc{num}(q,\delta _{\ell })$. The partial ordering on $%
Q_{\psi }$ is inclusion.

\begin{proposition}
\label{ccc}(i) For every $\delta \in E$ and $k\in \omega $, $D_{\delta
,k}=\{q\in Q_{\psi }:\delta \in \limfunc{cont}(q)$ and $k\leq \limfunc{num}%
(q,\delta )\}$ is dense in $Q_{\psi }$

(ii) $Q_{\psi }$ is c.c.c.
\end{proposition}

Before proving Proposition \ref{ccc}, we prove a lemma:

\begin{lemma}
\label{above}\bigskip Given $\{\delta _{0},...,\delta _{m}\}\in E$, integers 
$r_{\ell }^{\prime }$ for $\ell \leq m$ and a finite subset $I^{\prime }$ of 
$\omega _{1}$, there are integers $r_{\ell }^{\prime \prime }\geq r_{\ell
}^{\prime }$ for all $\ell \leq m$ and a finite subset $I^{\prime \prime }$
of $\omega _{1}$ containing $I^{\prime }$ such that for all $\ell \leq m$:

\begin{quote}
(a) $\eta _{\delta _{\ell }}(n)\in I^{\prime \prime }\Longleftrightarrow
n<r_{\ell }^{\prime \prime }$; and
\end{quote}

\begin{quote}
(b) for all $n<r_{\ell }^{\prime \prime }$, $u(\delta _{\ell },n)\in
\left\langle x_{\nu }:\nu \in I^{\prime \prime }\right\rangle $.
\end{quote}
\end{lemma}

\begin{proof}
The proof is by induction on $m\geq 0$. If $m=0$ we can take 
\[
r_{0}^{\prime \prime }=\max \{r_{0}^{\prime },\max \{k+1:\eta _{\delta
_{0}}(k)\in I^{\prime }\}\} 
\]
and take $I^{\prime \prime }$ to be a minimal extension of $I^{\prime }\cup
\{\eta _{\delta _{0}}(n):n<r_{0}^{\prime \prime }\}$ satisfying (b); then
(a) holds because $u(\delta ,n)\in \left\langle x_{\nu }:\nu <\eta _{\delta
}(n)\right\rangle $. If $m>0$, without loss of generality we can assume that 
$\delta _{0}<\delta _{1}<...<\delta _{m}$. Let 
\[
r_{m}^{\prime \prime }=\max \{r_{m}^{\prime },\max \{k+1:\eta _{\delta
_{m}}(k)\in I^{\prime }\},\min \{k:\eta _{\delta _{m}}(k)>\delta _{m-1}\}\}%
\text{.} 
\]
As in the case $m=0$, there exists $\tilde{I}$ containing $I^{\prime }$ such
that (a) and (b) hold for $\tilde{I}$ for $\ell =m$. Then apply the
inductive hypothesis to $\{\delta _{0},...,\delta _{m-1}\}$, $\tilde{I}$,
and the $r_{\ell }^{\prime }$ ($\ell <m$) to obtain $r_{\ell }^{\prime
\prime }$ for $\ell <m$ and a minimal $I^{\prime \prime }$.
\end{proof}

\smallskip \ For $q\in Q_{\psi }$ and $\alpha \in \omega _{1}$, let $%
q\upharpoonright \alpha $ denote the restriction of $q$ to 
\[
\limfunc{dom}(q)\cap (\{z_{\delta ,n}:\delta <\alpha ,n\in \omega \}\cup
\{x_{\nu }:\nu <\alpha \}). 
\]
Say that $\tau $ occurs in $q$ if $\tau \in \limfunc{cont}(q)$ or $x_{\tau
}\in \limfunc{dom}(q)$.

\smallskip\ 

\textsc{proof of proposition \ref{ccc}. }(i) Given $\delta \in E$, $k\in
\omega $ and $p\in Q_{\psi }$, we need ${q}${$\geq $}${p}$ such that $q\in
D_{\delta ,k}$. Let $\limfunc{cont}(p)=\{\delta _{0},...,\delta _{m}\}$. We
consider two cases. The first is that $\delta \in \limfunc{cont}(p)$, that
is, $\delta =\delta _{j}$ for some $j\leq m$. We can assume that $k>\limfunc{%
num}(p,\delta _{j})$. Apply Lemma \ref{above} with $I^{\prime }=I_{p}$, $%
r_{j}^{\prime }=k$, and $r_{\ell }^{\prime }=\limfunc{num}(p,\delta _{\ell })
$ for $\ell \neq j$ to obtain $I^{\prime \prime }$ and $r_{\ell }^{\prime
\prime }$. Then we can define $q$ to be the extension of $p$ with $\limfunc{%
cont}(q)=\limfunc{cont}(p)$ and $\limfunc{num}(q,\delta _{\ell })=r_{\ell
}^{\prime \prime }$ and $I_{q}=I^{\prime \prime }$. Since (\ref{5}) and (\ref
{7}) hold, we can inductively define $q(x_{\eta _{\delta _{\ell }}(i)})$ and 
$q(z_{\delta _{\ell },i+1})$ for $r_{\ell }^{\prime }\leq i<r_{\ell
}^{\prime \prime }$ (setting $q(x_{\nu })=0$ for $\nu \in I_{q}\setminus
\cup \{\limfunc{rge}(\eta _{\delta _{\ell }}:\ell \leq m\}$ if not already
defined) so that (\ref{6}) holds. Note that (\ref{7}) continues to hold.

The second case is when $\delta \notin \limfunc{cont}(p)$. Let $\delta
_{m+1}=\delta $. Choose $r_{\ell }^{\prime }$ for $\ell \leq m+1$ so that $%
r_{\ell }^{\prime }\geq \limfunc{num}(p,\delta _{\ell })$ for $\ell \leq m$
and such that (\ref{7}) holds, that is, $\eta _{\delta _{j}}(n)\neq \eta
_{\delta _{\ell }}(i)$ for all $n\geq r_{j}^{\prime }$ and $i\in \omega $
for all $j\neq \ell \in \{0,...,m+1\}$. Apply Lemma \ref{above} to \{$\delta
_{0},...,\delta _{m+1}\}$, $I_{p}$, and the $r_{\ell }^{\prime }$ to obtain $%
r_{\ell }^{\prime \prime }$ for $\ell \leq m+1$ and $I^{\prime \prime }$.
Let $I_{q}=I^{\prime \prime }$ and $\limfunc{num}(q,\delta _{\ell })=r_{\ell
}^{\prime \prime }$. For $\ell \leq m$ define $q(x_{\eta _{\delta _{\ell
}}(i)})$, $q(z_{\delta _{\ell },i+1})$ and $u(\delta _{\ell },i)$ for $%
\limfunc{num}(p,\delta _{\ell })\leq i<r_{\ell }^{\prime \prime }$ by
induction on $i$ as in the first case. Define $q(z_{\delta ,r_{m+1}^{\prime
\prime }+1})=0$ and define $q(z_{\delta ,n})$ for $n\leq r_{m+1}^{\prime
\prime }$ by ``downward induction'', i.e. 
\[
q(z_{\delta ,n})=2^{g(\delta ,n)}q(z_{\delta ,n+1})-q(x_{\eta _{\delta
}(n)})-q(u(\delta ,n))-\psi (w_{\delta ,n})\text{.} 
\]
(Setting $q(x_{\tau })=0$ where not already defined, we can assume $%
q(x_{\eta _{\delta }(n)})$ and $q(u(\delta ,n))$ are defined.)

(ii) Consider an uncountable subset $\{q_{\nu }:\nu \in \omega _{1}\}$ of $%
Q_{\psi }$. By the $\Delta $-system lemma we can assume that $\{\limfunc{cont%
}(q_{\nu }):\nu \in \omega _{1}\}$ forms a $\Delta $-system, i.e., there is
a finite subset $\Delta $ of $E$ such that for all $\nu \neq \mu $, $%
\limfunc{cont}(q_{\nu })\cap \limfunc{cont}(q_{\mu })=\Delta $. By
renumbering an uncountable subset, we can assume that for all $\nu $, if $%
\delta \in \limfunc{cont}(q_{\nu })\diagdown \Delta $, then $\delta >\nu $.
Furthermore, by passing to a subset and using (i) we can assume that if $%
\delta \in \limfunc{cont}(q_{\nu })\diagdown \Delta $ and $\eta _{\delta
}(n)<\nu $, then $n<\limfunc{num}(q_{\nu },\delta )$. By Fodor's Lemma we
can assume that there exists $\gamma \geq \max \Delta $ such that for all $%
\nu $ and $n$, if $\delta \in \limfunc{cont}(q_{\nu })$ and $\eta _{\delta
}(n)<\nu $, then $\eta _{\delta }(n)<\gamma $ and moreover such that if $%
\tau \in I_{q_{\nu }}$ and $\tau <\nu $, then $\tau <\gamma $. We can also
assume that for all $\mu ,\nu $, $q_{\mu }\upharpoonright \mu =q_{\nu
}\upharpoonright \nu $. If we pick $\mu <\nu $ such that $\gamma <\mu $ and
whenever $\tau $ occurs in $q_{\mu }$, then $\tau <\nu $, then we will have
that $q_{\mu }\cup $ $q_{\nu }$ $\in Q_{\psi }$. Notice that (\ref{7}) will
be satisfied: if $\delta \in \limfunc{cont}(q_{\mu })\setminus \Delta $ and $%
\rho \in \limfunc{cont}(q_{\nu })\setminus \Delta $ and $k\geq \limfunc{num}%
(q_{\mu },\delta )$ and $m\geq \limfunc{num}(q_{\nu },\rho )$, then $\mu
\leq $ $\eta _{\delta }(k)<\nu \leq \eta _{\rho }(m)$; moreover, if $i\in
\omega $ and $\eta _{\rho }(i)<\nu $, then $\eta _{\rho }(i)<\gamma <\mu
\leq \eta _{\delta }(k)$. Similarly it follows that (\ref{5}) holds. $%
\blacksquare $\ 

\smallskip\ 

Now $P=\left\langle P_{i},\dot{Q}_{i}:0\leq i<\omega _{2}\right\rangle $ is
defined to be a finite support iteration of length $\omega _{2}$ so that for
every $i\geq 1$ $\Vdash _{P_{i}}\dot{Q}_{i}=Q_{\dot{\psi}_{i}}$ where $%
\Vdash _{P_{i}}\dot{\psi}_{i}$ is a homomorphism$:K\rightarrow \mathbb{Z}$ and
the enumeration of names $\{\dot{\psi}_{i}:1\leq i<\omega _{2}\}$ is chosen
so that if $G$ is $P$-generic and $\psi \in V[G]$ is a homomorphism$%
:K\rightarrow \mathbb{Z}$, then for some $i\geq 1$, $\dot{\psi}_{i}$ is a name
for $\psi $ in $V^{P_{i}}$. Then $P$ is c.c.c. and in $V[G]$ every
homomorphism from $K$ to $\mathbb{Z}$ extends to one from $F$ to $\mathbb{Z}$.
This means that $\limfunc{Ext}(A,\mathbb{Z})=0$, that is, $A$ is a Whitehead
group (see, for example, \cite[p. 8]{EM}). We claim moreover that:

\begin{theorem}
\label{free}In $V[G]$ $A^{*}$ (= $\limfunc{Hom}(A,\mathbb{Z}\mathbf{)}$) is
free of cardinality $\aleph _{2}$.
\end{theorem}

As a consequence we can conclude:

\begin{corollary}
\label{freecor}In $V[G]$ $A$ is strongly non-reflexive.
\end{corollary}

\begin{proof}
Since $A^{*}$ is isomorphic to $\mathbb{Z}^{(\aleph _{2})}$,  $A^{**}$ is
isomorphic to $\mathbb{Z}^{\aleph _{2}}$ and hence not isomorphic to $A$
because its cardinality is different. We remark also that  $A^{**}$ is not
slender, but $A$ is slender since it is a Whitehead group --- see 
\cite[Prop. XII.1.3, p. 345]{EM}).
\end{proof}

\smallskip\ 

The next three sections are devoted to a proof of Theorem \ref{free}. The
fact that $A^{*}$ has cardinality $2^{\aleph _{1}}$ is a consequence of a
result of Chase \cite[Thm. 5.6]{C1}; by standard arguments it can be seen
that in $V[G]$ $2^{\aleph _{0}}=2^{\aleph _{1}}=\aleph _{2}$. Let $G_{\nu
}=\{p\upharpoonright \nu :p\in G\}$, so that $G_{\nu }$ is $P_{\nu }$%
-generic. To prove that $\limfunc{Hom}(A,\mathbb{Z}\mathbf{)}^{V[G]}$ is free,
it suffices to prove that:

(I) $\limfunc{Hom}(A,\mathbb{Z}\mathbf{)}^{V[G_{1}]}=0$;

(II) for every limit $\beta \leq \omega _{2}$, $\limfunc{Hom}(A,\mathbb{Z}%
\mathbf{)}^{V[G_{\beta }]}=\bigcup_{i<\beta }\limfunc{Hom}(A,\mathbb{Z}\mathbf{)%
}^{V[G_{i}]}$;

\noindent \noindent and

(III) for all $i<\omega _{2}$, $\limfunc{Hom}(A,\mathbb{Z}\mathbf{)}%
^{V[G_{i+1}]}/\limfunc{Hom}(A,\mathbb{Z}\mathbf{)}^{V[G_{i}]}$ is free, and in
fact is either $0$ or $\mathbb{Z}$.

We shall prove (I) immediately, and then prove the other two parts in the
next three sections.

\medskip\ 

\textsc{proof of (I):} Notice first that, by (\ref{1.2}), if $h\in \limfunc{%
Hom}(A,\mathbb{Z}\mathbf{)}$ and $h(x_{\mu })=0$ for all $\mu \in \omega _{1}$,
then $h$ is identically zero. So suppose, to obtain a contradiction, that
there exists a $Q_{0}$-name $\dot{h}$ and $p\in G_{1}$ such that 
\[
p\Vdash \dot{h}\in \limfunc{Hom}(A,\mathbb{Z}\mathbf{)}\wedge \dot{h}(x_{\mu
})=m 
\]
for some $\mu \in \omega _{1}$ and some non-zero integer $m$. Choose $d$
such that $2^{d}$ does not divide $m$. For each $\delta \in E$ there exists $%
p_{\delta }\geq p$ and $c_{\delta }\in \mathbb{Z}$ such that 
\[
p_{\delta }\Vdash \dot{h}(z_{\delta ,0})=c_{\delta }\text{.} 
\]
By Fodor's Lemma and a $\Delta $-system argument, there exist $\delta
_{1}\neq \delta _{2}>\mu $ such that $c_{\delta _{1}}=c_{\delta _{2}}$, and
if (for convenience of notation) we let $p^{i}=p_{\delta _{i}}$, $r_{\delta
_{1}}^{p_{{}}^{1}}=r_{\delta _{2}}^{p^{2}}=r$, $\eta _{\delta _{1}}(n)=\eta
_{\delta _{2}}(n)$, $u(\delta _{1},n)=u(\delta _{2},n)$ for all $n<r$ and $%
p^{1}$ and $p^{2}$ are compatible. Then there is a condition $q\in Q_{0}$
such that $p^{i}\leq q$ for $i=1,2$ and 
\[
q\Vdash \eta _{\delta _{1}}(r)=\eta _{\delta _{2}}(r)\wedge g(\delta
_{1},r)=d=g(\delta _{2},r)\wedge u(\delta _{1},r)=x_{\mu }\wedge u(\delta
_{2},r)=0\text{.} 
\]
Now consider a generic extension $V[G_{1}^{\prime }]$ where $q\in
G_{1}^{\prime }$. By subtracting (\ref{1.2}) for $n=r$ and $\delta =\delta
_{2}$ from (\ref{1.2}) for $n=r$ and $\delta =\delta _{1}$ and applying $h$
we obtain that (in $V[G_{1}^{\prime }]$) $2^{d}$ divides $h(u(\delta
_{1},r))-h(u(\delta _{2},r))=h(x_{\mu })-h(0)=m$. But this is a
contradiction of the choice of $d$. $\blacksquare $

\section{Preliminaries}

Before beginning the proof proper of (II) and (III), we prove a crucial
Proposition that we will need. For a fixed $m\in \omega $ and $S\subseteq E$%
, let $Z_{m}[S]$ denote the pure closure in $A$ of the subgroup generated by 
$\{z_{\delta ,m}+K:\delta \in S\}$. For $t\in \omega $, let $Z_{m,t}[S]$
denote $Z_{m}[S]+2^{t}A$.

\begin{proposition}
\label{zmt} In $V[G]$, for all $m,t\in \omega $ and all stationary $%
S\subseteq E$, $A/Z_{m,t}[S]$ is a finite group.
\end{proposition}

\begin{proof}
The proof is by contradiction. Suppose that $q^{*}\in G$ such that for some $%
m,t\in \omega $ and some $\dot{S}$%
\[
q^{*}%
\forces%
_{P}\dot{S}\text{ is a stationary subset of }E\text{ and }A/Z_{m,t}[\dot{S}]%
\text{ is infinite.}
\]
Let $S^{\prime }$ be the set of all $\delta \in E$ such that $q^{*}$ does
not force ``$\delta \notin \dot{S}$''; then $S^{\prime }\in V$ is a
stationary subset of $E$. For each $\delta \in S^{\prime }$ choose $%
p_{\delta }\geq q^{*}$ such that $p_{\delta }%
\forces%
\delta \in \dot{S}$.

We can assume that each $p_{\delta }$ satisfies:

\begin{quote}
($\dag $) $0\in \limfunc{dom}(p_{\delta })$; $\delta \in \limfunc{dom}%
(p_{\delta }(0))$; for each $j\in \limfunc{dom}(p_{\delta })$, $p_{\delta
}(j)$ is a function in $V$ and not just a name; $r_{\gamma }^{p_{\delta
}(0)} $ ($=$ $r_{\delta }$) is independent of $\gamma \in \limfunc{dom}%
(p_{\delta }(0))$; if $j\in \limfunc{dom}(p_{\delta })\setminus \{0\}$, $%
\gamma \in \limfunc{cont}(p_{\delta }(j))$ implies $\gamma \in \limfunc{dom}%
(p_{\delta }(0))$ and $\limfunc{num}(p_{\delta }(j),\gamma )$ ($=r_{\delta
,j}^{\prime } $) is $\leq r_{\delta }$ and independent of $\gamma $.
Moreover, if $\gamma \in \limfunc{dom}(p_{\delta }(0))$ and $\gamma >\delta $%
, then $\eta _{\gamma }(r_{\delta ,j}^{\prime })>\delta $.
\end{quote}

When we say that ``$\nu $ occurs in $p$'' we mean that $p(\nu )$ is
non-empty, or $\nu $ occurs in $p(j)$ for some $j>0$ or $\nu \in \limfunc{dom%
}(p(0))\cup \{\eta _{\gamma }(n):\gamma \in \limfunc{dom}(p(0))$, $%
n<r_{\gamma }^{p(0)}\}$, or $u(\gamma ,n)\notin \left\langle x_{\mu }:\mu
\neq \nu \right\rangle $ for some $n<r_{\gamma }^{p(0)}$. Without loss of
generality we can assume (passing to a subset of $S^{\prime }$) that

\begin{quote}
($\dag \dag $) there exists $\tau $ such that for all $\delta \in S^{\prime
} $, $\delta >\tau $ and every ordinal $<\delta $ which occurs in $p_{\delta
}$ is less than $\tau $; $\{\limfunc{dom}(p_{\delta }):\delta \in S^{\prime
}\}$ forms a $\Delta $-system, whose root we denote $C$ (i.e., $\limfunc{dom}%
(p_{\delta _{1}})\cap \limfunc{dom}(p_{\delta _{2}})=C=\{0,\mu _{1},...,\mu
_{d}\}$ for all $\delta _{1}\neq \delta _{2}$ in $S^{\prime }$); $r_{\delta
} $ ($=r^{*}$) and $r_{\delta ,j}^{\prime }$ ($=r_{j}^{\prime }$) are
independent of $\delta $. Moreover, for every $j\in C$, $\{\limfunc{dom}%
(p_{\delta }(j)):\delta \in S^{\prime }\}$ forms a $\Delta $-system and for
all $\delta _{1}\neq \delta _{2}$ in $S^{\prime }$, $p_{\delta _{1}}(j)$ and 
$p_{\delta _{2}}(j)$ agree on $\limfunc{dom}(p_{\delta _{1}}(j))\cap 
\limfunc{dom}(p_{\delta _{2}}(j))$, so $p_{\delta _{1}}(j)\upharpoonright
\delta _{1}=p_{\delta _{2}}(j)\upharpoonright \delta _{2}$. Also, $\limfunc{%
dom}(p_{\delta }(0))\cap \delta $ and $p_{\delta }(0)\upharpoonright (%
\limfunc{dom}(p_{\delta }(0))\cap \delta )$ are independent of $\delta $.
Finally, $\eta _{\delta }^{p_{\delta }(0)}(n)$ (=$\zeta _{n}$), $%
g^{p_{\delta }(0)}(\delta ,n)$ (= $g(n)$) and $u^{p_{\delta }(0)}(\delta ,n)$
are independent of $\delta $ for each $n<r^{*}$.
\end{quote}

Let $p^{*}$ denote the ``heart'' of $\{p_{\delta }:\delta \in S^{\prime }\}$%
; that is, $\limfunc{dom}(p^{*})=C$ and for all $\mu \in C$, $\limfunc{dom}%
(p^{*}(\mu ))=\limfunc{dom}(p_{\delta _{1}}(\mu ))\cap \limfunc{dom}%
(p_{\delta _{2}}(\mu ))$ (= $C_{\mu }$, say) for $\delta _{1}\neq \delta
_{2}\in S^{\prime }$; and $p^{*}(\mu )=p_{\delta _{1}}(\mu )\upharpoonright
C_{\mu }=p_{\delta _{1}}(\mu )\upharpoonright \delta _{1}$.

The conditions in ($\dag \dag $) insure that if $\delta _{1}<\delta _{2}$
are members of $S^{\prime }$ such that every ordinal which occurs in  $%
p_{\delta _{1}}$ is $<\delta _{2}$, then $p_{\delta _{1}}$ and $p_{\delta
_{2}}$ are almost compatible; however, there may be problems in determining
a value for $p_{\delta _{\ell }}(j)(x_{\zeta _{n}})$ for $r_{j}^{\prime
}\leq n<r^{*}$ (independent of $\ell =1,2$); it is because of these that the
following argument is necessary.

We can assume that $r^{*}$ $\geq m$ and that for all $\delta \in S^{\prime }$
$\delta \notin \limfunc{dom}(p^{*}(0))$. Choose  $M\geq t$ such that $%
g(n)\leq M$ for all  $n<r^{*}$. Let 
\[
N=2^{(r^{*}+1)M}\text{.}
\]

To obtain a contradiction, it suffices to prove that $p^{*}$ forces:

\begin{quote}
($\nabla $) $A/Z_{m,t}[S]$ is a group of cardinality $\leq $ $N^{d}$
\end{quote}

This is a contradiction since $p^{*}\geq q^{*}$. If $p^{*}$ does not force ($%
\nabla $), then there is a finite subset $\Theta $ of $\omega _{1}$ and a
condition $p^{**}\geq p^{*}$ such that $p^{**}$ forces

\begin{quote}
($\nabla \nabla $)$(\left\langle x_{\nu }:\nu \in \Theta \right\rangle
+Z_{m,t}[S])/Z_{m,t}[S]$ has cardinality $>N^{d}$.
 \end{quote}
(Note that it follows from (\ref{1.2}) that $A/p^{t}A$ is generated by $%
\left\langle x_{\nu }:\nu \in \Theta \right\rangle $.) We can assume that if 
$\nu $ occurs in $p^{**}$, then $\nu \in \Theta $. Let $T$ be the subset of $%
\left\langle x_{\nu }:\nu \in \Theta \right\rangle $ composed of all
elements of the form $\sum_{\nu \in \Theta }c_{\nu }x_{\nu }+x_{\beta }$
where $0\leq $ $c_{\nu }<2^{t}$. Let $\theta =$ $2^{|\Theta |t}$; so $T$ has 
$\theta >N^{d}$ elements; list them as $\{\tau _{\ell }:\ell <\theta \}$.
Now choose elements $\{\delta _{\ell }:\ell <\theta \}$ of $S^{\prime }$
listed in increasing order and such that the smallest, $\delta _{0}$, is
larger than $\max \Theta $ and such that every ordinal which occurs in $%
p_{\delta _{\ell }}$ is $<\delta _{\ell +1}$. Moreover, we can choose them
so that for any $\ell <\theta $, the ``common part'' of $p_{\delta _{\ell }}$
and $p^{**}$ is $p^{*}$; that is, $\limfunc{dom}(p_{\delta _{\ell }})\cap 
\limfunc{dom}(p^{**})=C=\limfunc{dom}(p^{*})$ and for all $\mu \in C$, $%
\limfunc{dom}(p_{\delta _{\ell }}(\mu ))\cap \limfunc{dom}(p^{**}(\mu ))=%
\limfunc{dom}(p^{*}(\mu ))$. (So, in particular, $\delta _{\ell }\notin 
\limfunc{dom}(p^{**}(0))$.)

Choose new ordinals $\alpha _{\ell }$ for $-1\leq \ell <\theta $ such that 
\[
\alpha _{-1}<\alpha _{0}<\delta _{0}<\alpha _{1}<\delta _{1}<...<\delta
_{\ell }<\alpha _{\ell +1}<\delta _{\ell +1}<... 
\]
Moreover, we make the choice so that for all $\ell $, $\alpha _{\ell }$ is
larger than any ordinal $<\delta _{\ell }$ which occurs in any $p_{\delta
_{k}}$. There is a condition $q_{0}\in Q_{0}$ which extends $p^{**}(0)$ and
each $p_{\delta _{\ell }}(0)$ ($\ell <\theta $) such that $q_{0}$ forces for
all $\ell <\theta $: 
\[
\eta _{\delta _{\ell }}(r^{*})=\alpha _{-1};\eta _{\delta _{\ell
}}(r^{*}+1)=\alpha _{\ell };u(\delta _{\ell },r^{*})=\tau _{\ell };\text{
and }g(\delta _{\ell },r^{*})=t. 
\]

\noindent This $q_{0}$ will force versions of (\ref{5}) and (\ref{7}) for
the $\delta _{\ell }$. Also choose $q_{0}$ to force values for $\eta
_{\gamma }(r^{*})$ \smallskip for any $\gamma \in \bigcup_{\ell <\theta }(%
\limfunc{dom}(p_{\delta _{\ell }}(0))-\{\delta _{\ell }\})$ so that (\ref{5}%
) and (\ref{7}) hold for $\allowbreak p_{\delta _{\ell }}(\mu )\cup
p_{\delta _{\ell ^{\prime }}}(\mu )$ for any $\ell ,\ell ^{\prime }<\theta $
and any $\mu \in C$. \ We claim that

\begin{quote}
(IV.1) There is a subset $W$ of $\{0,...,\theta -1\}$ of size at least $%
\theta \cdot N^{-d}$ and a condition $q\in P_{\omega _{2}}$ which extends $%
p^{**}$ and $p_{\delta _{\ell }}$ for every $\ell \in W$ and satisfies $%
q(0)=q_{0}$
\end{quote}

\noindent Assuming (IV.1), let us deduce a contradiction, which will prove
that $p^{*}$ forces ($\bigtriangledown $). Work in a generic extension $%
V[G^{\prime }]$ such that $q\in G^{\prime }$. For $\ell _{1}\neq \ell _{2}$
in $W$ we have $\tau _{\ell _{1}}-\tau _{\ell _{2}}\in \left\langle x_{\nu
}:\nu \in \Theta \right\rangle \cap Z_{m,t}[S]$ because 
\[
z_{\delta _{\ell _{1},r^{*}}}-z_{\delta _{\ell _{2}},r^{*}}=\tau _{\ell
_{1}}-\tau _{\ell _{2}}+2^{t}a
\]
for some $a\in A$ and, letting $e=\sum_{n=m}^{r^{*}-1}g(n)$, 
\[
2^{e}(z_{\delta _{\ell _{1},r^{*}}}-z_{\delta _{\ell _{2}},r^{*}})=z_{\delta
_{\ell _{1},m}}-z_{\delta _{\ell _{2}},m}\in Z_{m}[S]
\]
so since $Z_{m}[S]$ is pure-closed, $\tau _{\ell _{1}}-\tau _{\ell
_{2}}+2^{t}a\in Z_{m}[S]$, and hence $\tau _{\ell _{1}}-\tau _{\ell _{2}}\in
Z_{m,t}[S]$. Therefore $(\left\langle x_{\nu }:\nu \in \Theta \right\rangle
+Z_{m,t}[S])/Z_{m,t}[S]$ has cardinality at most 
\[
2^{|\Theta |t}/|W|=\theta /|W|\leq N^{d}
\]
which is a contradiction of the choice of $p^{**}$.

In order to prove (IV.1) we define  inductively, for $1\leq n\leq d+1$, a
condition $q_{n}\in P_{\mu _{n}}$ (where $\mu _{n}$ is as in the enumeration
of $C$ for $n\leq d$ and $\mu _{d+1}=\omega _{2}$) such that $q_{n}\geq
p^{**}\upharpoonright \mu _{n}$ and for $n^{\prime }<n$, $%
q_{n}\upharpoonright \mu _{n^{\prime }}\geq q_{n^{\prime }}$. We also define
a subset $W_{n}$ of $W$ of size at least $\theta \cdot N^{-(n-1)}$ such that
for each $\ell \in W_{n}$, $q_{n}\geq p_{\delta _{\ell }}\upharpoonright \mu
_{n}$. (So in the end we let $q=q_{d+1}$ and $W=W_{d+1}$.)

To begin, let $W_{1}=\{0,...,\theta -1\}$ and let $q_{1}$ be any common
extension of $q_{0}$ and the $p_{\delta _{\ell }}\upharpoonright \mu _{1}$.
(There is no problem finding such an extension.) Suppose now that $q_{n}$
and $W_{n}$ have been defined for some $n\geq 1$. Choose $\tilde{q}_{n}\geq
q_{n}$ in $P_{\mu _{n}}$ such that $\tilde{q}_{n}$ decides for all $\ell \in
W_{n}$ the value of $\psi _{\mu _{n}}(w_{\gamma ,k})$ for all $\gamma \in 
\limfunc{dom}(p_{\delta _{\ell }(0)})$ and all $k\leq r^{*}$. For each $\ell
\in W_{n}$ fix $s_{\ell }\in V$ such that $s_{\ell }\in Q_{\mu _{n}}$
extends $p_{\delta _{\ell }}(\mu _{n})$ and satisfies $\limfunc{num}(s_{\ell
},\delta _{\ell })=r^{*}+1$, $s_{\ell }(x_{\alpha _{-1}})<2^{M}$, and $%
s_{\ell }(x_{\zeta _{k}})<2^{M}$ for all $k<r^{*}$. (This is possible by the
proof of Proposition \ref{ccc} since we only need to find solutions to the
equations (\ref{1.15}) modulo $2^{M}$since $g^{p_{\delta _{\ell
}}(0)}(\delta _{\ell },k)\leq M$ for $k\leq r^{*}$.)

Define an equivalence relation $\equiv _{n}$ on $W_{n}$ by: $\ell _{1}\equiv
_{n}\ell _{2}$ iff $s_{\ell _{1}}\cup s_{\ell _{2}}$ is a function. By
choice of $M$ and $N$, there is an equivalence class, $W_{n+1}$, of size at
least $|W_{n}|/N$. For $\ell \in W_{n+1}$ we can define a common extension $%
q_{n+1}(\mu _{n})$ of $p^{**}(\mu _{n})$ and the $p_{\delta _{\ell }}(\mu
_{n})$ and let $q_{n+1}\upharpoonright \mu _{n}=\tilde{q}_{n}$. This
completes the inductive construction.
\end{proof}

For any abelian group $H$, let $\nu (H)$ be the Chase radical of $H$: the
intersection of the kernels of all homomorphisms of $H$ into an $\aleph _{1}$%
-free group (cf. \cite{C2}). Then $H/\nu (H)$ is $\aleph _{1}$-free (%
\cite[Prop. 1.2]{C2}, \cite[p. 290]{EM}). Let $\limfunc{cl}(Z_{m}[S])$ be
defined by: $\limfunc{cl}(Z_{m}[S])/Z_{m}[S]=\nu (A/Z_{m}[S])$, so in
particular $A/\limfunc{cl}(Z_{m}[S])$ is $\aleph _{1}$-free. Notice also
that every homomorphism from $A$  to $\mathbb{Z}$ is determined on $\limfunc{cl}%
(Z_{m}[S])$ by its values on $\{z_{\delta ,m}+K:\delta \in S\}$.

\begin{corollary}
\label{zmt1} In $V[G]$, for all $m\in \omega $ and stationary $S\subseteq E$%
, $A/\limfunc{cl}(Z_{m}[S])$ is a finite rank free group.
\end{corollary}

\begin{proof}
If not, then since $A/\limfunc{cl}(Z_{m}[S])$ is $\aleph _{1}$-free, it
contains a free pure subgroup of countably infinite rank. Let $\{a_{n}+%
\limfunc{cl}(Z_{m}[S]):n\in \omega \}$ be a basis of such a subgroup. For
any $n\neq m$, $a_{n}+Z_{m,1}[S]\neq a_{m}+Z_{m,1}[S]$ since $2$ does not
divide $a_{n}-a_{m}$ mod $Z_{m}[S]$ (or even mod $\limfunc{cl}(Z_{m}[S])$).
Therefore $\{a_{n}+$ $Z_{m,1}[S]:n\in \omega \}$ is an infinite subset of $%
A/Z_{m,1}[S]$, which contradicts Proposition \ref{zmt}.
\end{proof}

For the next Corollary we will need the following:

\begin{lemma}
\label{abs}The Chase radical, $\nu (H)$, of a torsion-free group $H$ is
absolute for generic extensions.
\end{lemma}

\begin{proof}
We give an absolute construction of $\nu (H)$ using the fact that a
torsion-free group is $\aleph _{1}$-free if and only if every finite rank
subgroup is finitely-generated (cf. \cite[Thm. 19.1]{F}), that is, if and
only if the pure closure of every finitely-generated subgroup is
finitely-generated. For any group $H^{\prime }$, let $\mu (H^{\prime })$ be
the sum of all finite rank subgroups $G$ of $H^{\prime }$ which are not free
but are such that every subgroup of $G$ of smaller rank is free; it is easy
to see that for such $G$, $\nu (G)=G$ and hence $\mu (H^{\prime })\subseteq
\nu (H^{\prime })$. Moreover, the definition of $\mu (H^{\prime })$ is
absolute. Now define $\nu _{\beta }(H)$ by induction: $\nu _{0}(H)=0$, $\nu
_{\alpha +1}(H)/\nu _{\alpha }(H)=\mu (H/\nu _{\alpha }(H)),$and for limit
ordinals $\beta $, $\nu _{\beta }(H)=\cup _{\alpha <\beta }\nu _{\alpha }(H)$%
. It follows by induction that $\nu _{\alpha }(H)\subseteq \nu (H)$ for all $%
\alpha \leq \omega _{1}$. We claim that $\nu (H)=\nu _{\omega _{1}}(H)$; it
suffices to prove that $H/\nu _{\omega _{1}}(H)$ is $\aleph _{1}$-free. If
not, then there is a finite rank subgroup of $H/\nu _{\omega _{1}}(H)$ which
is not finitely-generated. We can choose one, $G$, of minimal rank, so all
of its subgroups of smaller rank are free; say $G$ is the pure closure of $%
\{a_{1}+\nu _{\omega _{1}}(H),...,a_{n}+\nu _{\omega _{1}}(H)\}$; but then
for some $\alpha <\omega _{1}$, the pure closure of $\{a_{1}+\nu _{\alpha
}(H),...,a_{n}+\nu _{\alpha }(H)\}$ is not free, but still has the property
that every subgroup of smaller rank is free; hence $\{a_{1},...,a_{n}\}%
\subseteq \nu _{\alpha +1}(H)$, which is a contradiction.
\end{proof}

\begin{corollary}
\label{staty}If $h\in \limfunc{Hom}(A,\mathbb{Z}\mathbf{)}^{V[G]}$ and for some 
$i\in \omega _{2}$, $m\in \omega $ and some stationary set $S\in V[G_{i}]$,
the sequence $(h(z_{\delta ,m}+K):\delta \in S)$ belongs to $V[G_{i}]$, then 
$h$ belongs to $V[G_{i}]$.
\end{corollary}

\begin{proof}
Suppose $h$, $S$ and $m$ are as in the hypotheses. First we claim that $%
h\upharpoonright Z_{m}[S]$ belongs to $V[G_{i}]$. Indeed we can define $%
h\upharpoonright Z_{m}[S]$ in $V[G_{i}]$ as follows: $h(a)=k$ if for some $%
n\neq 0$, $na$ belongs to the subgroup generated by $\{z_{\delta
,m}+K:\delta \in S\}$ and $h(na)=nk$; and otherwise $h(a)=\xi $ for some
fixed $\xi \notin \mathbb{Z}$. In fact, the second case never occurs because $%
h(a)$ is defined in $V[G]$. Next we claim that $h\upharpoonright \limfunc{cl}%
(Z_{m}[S])$ belongs to $V[G_{i}]$. The proof is similar in principle, using
the inductive construction of $\limfunc{cl}(Z_{m}[S])$ given by the proof of
Lemma \ref{abs}. But then by Corollary \ref{zmt1}, $h$ is determined by only
finitely many more values, so also $h$ belongs to $V[G_{i}]$.
\end{proof}

\smallskip \ 

For the next corollary we introduce some notation that will be used in
section 4. Let $\varphi _{i}\in $ $V[G_{i+1}]$ denote the generic function
for $Q_{i}$; thus $\varphi _{i}$ is a homomorphism: $F\rightarrow {\mathbb Z}
$ extending $\psi _{i}:K\rightarrow {\mathbb Z}$, where $\psi _{i}$ is the
interpretation in $V[G_{i}]$ of the name $\dot{\psi}_{i}$. The canonical map
$\operatorname{Hom}(K, {\mathbb Z} )\rightarrow $ $\operatorname{Ext}(A,{%
\mathbb Z})$ sends $\psi _{i}$ to a short exact sequence

\[
\mathcal{E}_{i}:\text{ }0\rightarrow {\mathbb Z}\stackrel{\iota }{%
\longrightarrow }B_{i}\stackrel{\pi }{\longrightarrow }A\rightarrow 0 
\]

\noindent and there is a commuting diagram 
\[
\begin{array}{cclllcccc}
0 & \rightarrow  & K & \stackrel{\iota ^{\prime }}{\longrightarrow } & F & 
\longrightarrow  & A & \rightarrow  & 0 \\ 
&  & \downarrow _{{}}^{\psi _{i}} &  & \downarrow ^{\sigma _{i}} &  & 
\downarrow ^{1_{A}} &  &  \\ 
0 & \rightarrow  & {\mathbb Z} & \stackrel{\iota }{\longrightarrow } & B_{i}
& \stackrel{\pi }{\longrightarrow } & A & \rightarrow  & 0
\end{array}
\]
where $\iota $ and $\iota ^{\prime }$ are inclusion maps. Moreover, for all $%
z\in F$, ($\pi \circ \sigma _{i})(z)=z+K\in A$. Then $\varphi _{i}$ gives
rise to a splitting $\rho _{i}\in \limfunc{Hom}(B_{i},{\mathbb Z}%
)^{V[G_{i+1}]}$ defined by $\rho _{i}(\sigma _{i}(z))=\varphi _{i}(z)$. Thus 
$\rho _{i}\circ \iota =1_{\mathbb{Z}}$ (the identity on ${\mathbb Z}$) and also 
$\rho _{i}\circ \sigma _{i}=\varphi _{i}$.


\begin{corollary}
\label{Bi} If $f\in \limfunc{Hom}(B_{i},\mathbb{Z}\mathbf{)}^{V[G]}$ and for
some $m\in \omega $ and some stationary set $S\in V[G_{i}]$, the sequence $%
(f(\sigma _{i}(z_{\delta ,m})):\delta \in S)$ belongs to $V[G_{i}]$, then $f$
belongs to $V[G_{i}]$.
\end{corollary}

\begin{proof}
If $Z^{\prime }$ is defined to be the pure subgroup of $B_{i}$ generated by $%
\{\sigma _{i}(z_{\delta ,m}):\delta \in S\}\cup \{\iota (1)\}$ and $%
C^{\prime }$ is such that $\nu (B_{i}/Z^{\prime })=C^{\prime }/Z^{\prime }$,
then $\pi $ induces an isomorphism of $C^{\prime }/\limfunc{rge}(\iota )$
with cl($Z_{m}[S])$. Hence $B_{i}/C^{\prime }\cong A/$ cl($Z_{m}[S])$ is
finite rank free; therefore, arguing as in Corollary \ref{staty},  $f$
belongs to $V[G_{i}]$.
\end{proof}

\section{Proof of (II)}

We divide the proof of (II) into three cases according to the cofinality of $%
\beta $. The case of cofinality $\omega _{2}$ (i.e., $\beta =\omega _{2}$)
is trivial since any function from $A$ (which has cardinality $\aleph _{1}$)
to $\mathbb{Z}$ must belong to $V[G_{i}]$ for some $i<\beta $.

Let $\dot{h}$ be a $P_{\beta }$-name and $p\in G_{\beta }$ such that 
\[
p\Vdash _{P_{\beta }}\dot{h}\in \limfunc{Hom}(A,\mathbb{Z}\mathbf{)}\text{.} 
\]
Then for each $\delta \in E$ there is $p_{\delta }\in G_{\beta }$ and $%
k_{\delta }\in \mathbb{Z}$ such that $p_{\delta }\geq p$ and $p_{\delta }\Vdash
_{P_{\beta }}\dot{h}(z_{\delta ,0}+K)=k_{\delta }$.

Suppose that the cofinality of $\beta $ is $\omega $, and fix an increasing
sequence $(\beta _{n}:n\in \omega )$ whose sup is $\beta $. Then there is $%
n\in \omega $ and a stationary subset $S_{1}$ of $E$, belonging to $%
V[G_{\beta }]$, such that for $\delta \in S_{1}$, $p_{\delta }\in G_{\beta
_{n}}$. Without loss of generality there is $p^{*}\in G_{\beta n}$ such that 
$p^{*}$ forces 
\[
S=\{\delta \in E:\exists p_{\delta }\in \dot{G}_{\beta _{n}}\text{ and }%
k_{\delta }\in \mathbb{Z}\text{ s.t. }p_{\delta }\Vdash _{P_{\beta }}\dot{h}%
(z_{\delta ,0}+K)=k_{\delta }\}\text{ is stationary}.
\]
Then $(h(z_{\delta ,0}+K):\delta \in S)$ belongs to $V[G_{\beta _{n}}]$, so $%
h\in \limfunc{Hom}(A,\mathbb{Z}\mathbf{)}^{V[G_{\beta _{n}}]}$, by Corollary 
\ref{staty}.

\smallskip\ 

The final, and hardest, case is when the cofinality of $\beta $ is $\omega
_{1}$. Fix an increasing continuous sequence $(\beta _{\nu }:\nu <\omega
_{1})$ whose sup is $\beta $. Then there is $\nu \in \omega _{1}$ and a
stationary subset $S$ of $\omega _{1}$ such that for $\delta \in S$, $%
p_{\delta }\upharpoonright \beta _{\delta }\in G_{\beta _{\nu }}$.

For any $t\geq 1$, $(Z_{0,t}[S\setminus \alpha ]:\alpha <\omega _{1})$ is a
non-increasing sequence of groups. Since the groups $A/Z_{0,t}[S\backslash
\alpha ]$ are finite, it follows that there is a countable ordinal $\alpha
_{t}$ such that for $\gamma ,\alpha \geq \alpha _{t}$, $Z_{0,t}[S\setminus
\alpha ]=Z_{0,t}[S\setminus \gamma ]$. Therefore there is a countable
ordinal $\alpha _{*}$ and a countable subset $Y$ of $A$ such that for all $%
t\geq 1$, $\alpha _{t}\leq \alpha _{*}$, and $Y$ contains representatives of
all the elements of $A/Z_{0,t}[S\setminus \alpha _{*}]$. Increasing $\nu $
if necessary, we can assume that we can compute $h(y)$ in $V[G_{\beta _{\nu
}}]$ for all $y\in Y$.

We claim that $h$ belongs to $V[G_{\beta _{\nu }}]$. In pursuit of a
contradiction, suppose that there are $a\in A$, conditions $q_{1},q_{2}\in
P_{\beta }/G_{\beta _{\nu }}$ and integers $c_{1}\neq c_{2}$ such that $%
q_{\ell }\Vdash _{P_{\beta }}\dot{h}(a)=c_{\ell }$ for $\ell =1,2$. Choose $%
t $ sufficiently large such that $2^{t}$ does not divide $c_{2}-c_{1}$ and
choose $\mu <\omega _{1}$ such that $q_{1},q_{2}\in P_{\beta _{\mu }}$. For
some $y\in Y$, $a-y\in Z_{0,t}[S\setminus \alpha _{*}]=Z_{0,t}[S\setminus
\beta _{\mu }]$. Thus 
\[
a-y=z+2^{t}a^{\prime } 
\]
for some $a^{\prime }\in A$ and $z$ in the pure closure of the subgroup
generated by $\{z_{\delta _{j},0}:j=1,...,n\}$ for some $\delta
_{1},...,\delta _{n}\in S\setminus \beta _{\mu }$. For $\ell =1,2$ there is
an upper-bound $r_{\ell }\in P_{\beta }$ of $\{q_{\ell },p_{\delta
_{1}},...,p_{\delta _{n}}\}$. Then $r_{1}$ and $r_{2}$ force the same value, 
$b$, to $h(z)$ (because they are both $\geq p_{\delta _{j}}$ for $j=1,...,n$%
) and the same value, $k$, to $h(y)$ (because it is determined in $%
V[G_{\beta _{\nu }}]$). Therefore 
\[
r_{\ell }\Vdash 2^{t}\text{ divides }c_{\ell }-k-b\text{.} 
\]
So for $\ell =1,2$, the integer $c_{\ell }-k-b$ is divisible by $2^{t}$. But
this contradicts the choice of $t$.

\section{Proof of (III)}

We continue with the notation from the end of section 2; so $\mathcal{E}%
_{i}\in $ $V[G_{i}]$. Suppose that $\mathcal{E}_{i}$ represents a torsion
element of $\limfunc{Ext}(A,\mathbb{Z})$, of order $e\geq 1$, that is, there is
a homomorphism $g_{i}:B_{i}\rightarrow \mathbb{Z}$ such that $%
g_{i}\upharpoonright \mathbb{Z}=e1_{\mathbb{Z}}$, or more precisely, $g_{i}\circ
\iota =e1_{\mathbb{Z}}$. (We consider the zero element to be torsion of order $1
$.) Then $e\rho _{i}-g_{i}$ is a homomorphism from $B_{i}$ to $\mathbb{Z}$
which is identically $0$ on $\mathbb{Z}$, so it induces a homomorphism $\theta
_{i}\in \limfunc{Hom}(A,\mathbb{Z}\mathbf{)}^{V[G_{i+1}]}$ (that is, $\theta
_{i}\circ \pi =e\rho _{i}-g_{i}$) which is a new element of $A^{*}$ --- that
is, it is not in $V[G_{i}]$. To prove (III) it will suffice to prove that if
there is an element $h$ of $A^{*}$ which is in $V[G_{i+1}]$ but not in $%
V[G_{i}]$, then $\mathcal{E}_{i}$ is torsion, and in that case $h$ is an
integral multiple of $\theta _{i}$ modulo $(A^{*})^{V[G_{i}]}$.

Given such an $h$, let $h^{\prime }=h\circ \pi :B_{i}\rightarrow \mathbb{Z}$.
Clearly $h^{\prime }\in $ $V[G_{i+1}]-V[G_{i}]$. We claim that:

\begin{quote}
(III.1) For some integer $c$, $h^{\prime }-c\rho _i$ belongs to $V[G_i]$.
\end{quote}

Let us see first why this Claim implies the desired conclusion. Note that $%
c\neq 0$ since $h^{\prime }$ does not belong to $V[G_{i}]$. Since $%
(h^{\prime }-c\rho _{i})\upharpoonright \mathbb{Z}=-c1_{\mathbb{Z}}$, we conclude
that in $V[G_{i}]$, $\mathcal{E}_{i}$ is torsion, of order $e$ dividing $-c$%
; let $g_{i}\in \limfunc{Hom}(B_{i},\mathbb{Z}\mathbf{)}^{V[G_{i}]}$ such that $%
g_{i}\upharpoonright \mathbb{Z}=e1_{\mathbb{Z}}$. Let $\theta _{i}$ be induced by $%
e\rho _{i}-g_{i}$, as above. Say $c=ne$; then $h^{\prime }-c\rho _{i}+ng_{i}$
belongs to $V[G_{i}]$ and is identically $0$ on $\mathbb{Z}$ so it induces a
homomorphism $f\in \limfunc{Hom}(A,\mathbb{Z}\mathbf{)}^{V[G_{i}]}$. By
composing both sides with $\pi $ one sees that $h=n\theta _{i}+f$.

\smallskip\ 

We shall now work on the proof of (III.1). Let $F^{\prime }$ be the subgroup
of $F$ generated by $\{x_{\nu }:\nu <\omega _{1}\}$. We work in $V[G_{i}]$.
For any countable ordinal $\alpha \in \omega _{1}-E$, define 
\[
Q_{i,\alpha }=\{q\in Q_{i}:z_{\delta ,n}\in \limfunc{dom}(q)\Rightarrow
\delta <\alpha \text{ and }x_{\nu }\in \limfunc{dom}(q)\Rightarrow \nu
<\alpha \}. 
\]
Then $Q_{i,\alpha }$ is a complete subforcing of $Q_{i}$. In particular, 
\[
V[G_{i+1}]=V[G_{i}][G_{i+1,\alpha }][H_{i+1,\alpha }] 
\]
where $G_{i+1,\alpha }$ is $Q_{i,\alpha }$-generic over $V[G_{i}]$ and $%
H_{i+1,\alpha }$ is $Q_{i}/G_{i+1,\alpha }$-generic over $%
V[G_{i}][G_{i+1,\alpha }]$. We claim:

\begin{quote}
(III.2) There is a countable ordinal $\alpha \in \omega _{1}-E$ such that in 
$V[G_{i}][G_{i+1,\alpha }]$ there is an assignment to every $y\in F^{\prime
} $ of a function $\xi _{y}:\mathbb{Z}\rightarrow \mathbb{Z}$ such that for all $%
y\in F^{\prime }$ $%
\forces%
_{Q_{i}/G_{i+1,\alpha }}\dot{h}(y+K)=\xi _{y}(\dot{\varphi}_{i}(y))$.
\end{quote}

Let us see first why this implies (III.1). First we assert that the
following consequence of (III.2) holds in $V[G_i][G_{i+1,\alpha }]$:

\begin{quote}
(III.2.1) There is an integer $c$ such that for every $k\in \mathbb{Z}$, and
every $\beta \geq \alpha $, $\xi _{x_{\beta }}(k)-\xi _{x_{\beta }}(0)=kc$.
\end{quote}

To see this, let $\beta ,\gamma \geq \alpha $ with $\beta \neq \gamma $, and
let $k\in \mathbb{Z}$. By the proof of Proposition \ref{ccc}, there are
conditions $q_{1},q_{2}\in Q_{i}/G_{i+1,\alpha }$ such that 
\[
q_{1}%
\forces%
\dot{\varphi}_{i}(x_{\beta })=k\wedge \dot{\varphi}_{i}(x_{\gamma })=0
\]
and 
\[
q_{2}%
\forces%
\dot{\varphi}_{i}(x_{\beta })=0\wedge \dot{\varphi}_{i}(x_{\gamma })=k\text{.%
}
\]
Let $y=x_{\beta }+x_{\gamma }$. By (III.2) and the fact that $\dot{h}$ and $%
\varphi _{i}$ are homomorphisms, 
\[
q_{1}%
\forces%
\xi _{x_{\beta }}(k)+\xi _{x_{\gamma }}(0)=\xi _{y}(k)
\]
which implies that $\xi _{x_{\beta }}(k)+\xi _{x_{\gamma }}(0)=\xi _{y}(k)$
holds in $V[G_{i}][G_{i+1,\alpha }]$. Similarly, reasoning with $q_{2}$, we
can conclude that $\xi _{x_{\gamma }}(k)+\xi _{x_{\beta }}(0)=\xi _{y}(k)$
holds in $V[G_{i}][G_{i+1,\alpha }]$. Thus $\xi _{x_{\beta }}(k)-\xi
_{x_{\beta }}(0)=\xi _{x_{\gamma }}(k)-\xi _{x_{\gamma }}(0)$ in $%
V[G_{i}][G_{i+1,\alpha }]$; we denote this value by $\Xi (k)$. If we can
prove that for all $k$, $\Xi (k)=k\Xi (1)$, then we can let $c=\Xi (1)$.
Again,  let $\beta ,\gamma \geq \alpha $ with $\beta \neq \gamma $ and this
time let $y=kx_{\beta }+x_{\gamma }$. Using conditions $q_{3}%
\forces%
\dot{\varphi}_{i}(x_{\beta })=1\wedge \dot{\varphi}_{i}(x_{\gamma })=0$ and $%
q_{2}%
\forces%
\dot{\varphi}_{i}(x_{\beta })=0\wedge \dot{\varphi}_{i}(x_{\gamma })=k$ we
conclude that 
\[
k\xi _{x_{\beta }}(1)+\xi _{x_{\gamma }}(0)=k\xi _{x_{\beta }}(0)+\xi
_{x_{\gamma }}(k)
\]
from which it follows that $k\Xi (1)=\Xi (k)$. This proves (III.2.1)

\smallskip\ 

Now work in $V[G_{i+1}]$; we have 
\[
h(x_{\beta }+K)=\xi _{x_{\beta }}(\varphi _{i}(x_{\beta }))=c\varphi
_{i}(x_{\beta })+\xi _{x_{\beta }}(0)
\]
for $\beta \geq \alpha $. Since $((h^{\prime }-c\rho _{i})\circ \sigma
_{i})(x)=h(x+K)-c\varphi _{i}(x)$ for $x\in F^{\prime }$, it follows that $%
h^{\prime }-c\rho _{i}\upharpoonright \{\sigma _{i}(x_{\beta }):$ $\beta
\geq \alpha \}$ belongs to $V[G_{i}][G_{i+1,\alpha }]$. Moreover, for $\beta
<\alpha $, $\varphi _{i}(x_{\beta })$ is determined in $V[G_{i}][G_{i+1,%
\alpha }]$, and hence so are $h(x_{\beta }+K)=\xi _{x_{\beta }}(\varphi
_{i}(x_{\beta }))$ and  $(h^{\prime }-c\rho _{i})(\sigma _{i}(x_{\beta }))$.
Therefore $h^{\prime }-c\rho _{i}$ belongs to $V[G_{i}][G_{i+1,\alpha }]$
(since it is determined by its values on $\{\sigma _{i}(x_{\beta }):\beta
\in \omega _{1}\}\cup \{\iota (1)\}$). Let $f=h^{\prime }-c\rho _{i}$. For
each $\delta \in E$, there exist $p_{\delta }\in G_{i+1,\alpha }$ and $%
k_{\delta }\in \mathbb{Z}$ such that 
\[
p_{\delta }%
\forces%
_{Q_{i}}\dot{f}(\sigma _{i}(z_{\delta ,m}))=k_{\delta }.
\]
Since $Q_{i,\alpha }$ and $\mathbb{Z}$ are countable, there exist $p\in
G_{i+1,\alpha }$, $k\in \mathbb{Z}$, and a stationary $S\in V[G_{i}]$ such that
for $\delta \in S$, $p%
\forces%
_{Q_{i}}\dot{f}(\sigma _{i}(z_{\delta ,m}))=k$. Then the (constant) sequence 
$(f(\sigma _{i}(z_{\delta ,m})):\delta \in S)$ belongs to $V[G_{i}]$, so by
Corollary \ref{Bi}, $f$ belongs to $V[G_{i}]$.

\smallskip\ 

So it remains to prove (III.2). Work in $V[G_{i}]$. Let 
\[
D_{i,\alpha }=\{q\in Q_{i}:\forall \delta \in (\limfunc{cont}(q)-\alpha
)\forall n\in \omega [(\eta _{\delta }(n)<\alpha )\Rightarrow x_{\eta
_{\delta }(n)}\in \limfunc{dom}(q)]\}\text{.} 
\]
Then $D_{i,\alpha }$ is a dense subset of $Q_{i}$. We claim that it is true
in $V[G_{i}]$ that:

\begin{quote}
(III.3) there is a countable ordinal $\alpha \in \omega _{1}-E$ such that
for every $y\in F^{\prime }$, $t,c_{1},c_{2}\in \mathbb{Z}$, and $%
q_{1},q_{2}\in D_{i,\alpha }$ with $q_{1}\upharpoonright \alpha
=q_{2}\upharpoonright \alpha $, if 
\[
q_{\ell }%
\forces%
_{Q_{i}}\dot{\varphi}_{i}(y)=t\wedge \dot{h}(y+K)=c_{\ell } 
\]
for $\ell =1,2$, then $c_{1}=c_{2}$.
\end{quote}

Clearly this implies (III.2). Indeed, we define $\xi _{y}(t)$ to be $c$ if
there is a $q\in D_{i,\alpha }$ such that $q\upharpoonright \alpha \in
G_{i+1,\alpha }$ and $q%
\forces%
_{Q_{i}/G_{i+1,\alpha }}\dot{\varphi}_{i}(y)=t\wedge \dot{h}(y+K)=c$ and
otherwise $\xi _{y}(t)=0$. By (III.3), $\xi _{y}(t)$ is well-defined.

\smallskip\ 

\textsc{proof of (III.3).} The proof is by contradiction and uses some of
the methods of the proof of Proposition \ref{zmt}. So suppose that for every 
$\alpha \in \omega _{1}-E$ there are $y^{\alpha }\in F^{\prime }$, $%
t^{\alpha },c_{1}^{\alpha },c_{2}^{\alpha }\in \mathbb{Z}$, and $q_{1}^{\alpha
},q_{2}^{\alpha }\in D_{i,\alpha }$ such that $q_{1}^{\alpha
}\upharpoonright \alpha =q_{2}^{\alpha }\upharpoonright \alpha $ and $%
q_{\ell }^{\alpha }%
\forces%
_{Q_{i}}\dot{\varphi}_{i}(y^{\alpha })=t^{\alpha }\wedge \dot{h}^{\prime
\prime }(y)=c_{\ell }^{\alpha }$ where $c_{1}^{\alpha }\neq c_{2}^{\alpha }$
for $\ell =1,2$. Then, by Fodor's Lemma and counting, there is a $p_{0}\in
G_{i}$, $t,c_{1},c_{2}\in \mathbb{Z}$, $\tilde{q}$ $\in V$ and names $\dot{S}$, 
$\dot{q}_{\ell }^{\alpha }$, $\dot{y}^{\alpha }$ such that 
\[
\begin{array}{c}
p_{0}%
\forces%
_{P_{i}}\dot{S}\text{ is a stationary subset of }\omega _{1}-E\text{ s.t.
for all }\alpha \in \dot{S}, \\ 
\text{\ }\dot{t}^{\alpha }=t,\dot{c}_{1}^{\alpha }=c_{1},\dot{c}_{2}^{\alpha
}=c_{2}\mathbf{\ }\text{and }\dot{q}_{1}^{\alpha }\upharpoonright \alpha =%
\dot{q}_{2}^{\alpha }\upharpoonright \alpha =\tilde{q}
\end{array}
\]
and moreover such that $p_{0}$ forces the names to be a counterexample to
(III.3), as above.

There is a stationary subset $S^{\prime }\subseteq \omega _{1}-E$ such that
for every $\alpha \in S^{\prime }$ there is a condition $p_{\alpha }\geq
p_{0}$ in $P_{i}$ which forces $\alpha \in \dot{S}$ and forces values
(elements of $V$) to $\dot{q}_{\ell }^{\alpha }$ and to $\dot{y}^{\alpha }$.
Moreover, we can suppose that the $p_{\alpha }\cup \{(i,q_{\ell }^{\alpha
})\}\in P_{i+1}$ ($\ell =1,2$) are as in ($\dag $) [cf. proof of Proposition 
\ref{zmt}] and that $\{p_{\alpha }:\alpha \in S^{\prime }\}$ is as in ($\dag
\dag $) [with $\alpha $ in place of $\delta $, but since $\alpha \notin E$,
the last sentence does not apply]. Let $p^{*}$ be the heart of $\{p_{\alpha
}:\alpha \in S^{\prime }\}$. We can also assume that $\{q_{\ell }^{\alpha
}:\alpha \in S^{\prime }\}$ forms a $\Delta $-system with heart $\tilde{q}$
(for $\ell =1,2$).

For each $\delta \in E$, there is $\hat{p}_{\delta }\in P_{i+1}$ and $%
k_{\delta }\in \mathbb{Z}$ such that $\hat{p}_{\delta }\upharpoonright i\geq
p^{*}$, $\hat{p}_{\delta }(i)\geq \tilde{q}$ and $\hat{p}_{\delta }\Vdash 
\dot{h}(z_{\delta ,0})=k_{\delta }$. There is a stationary $\hat{S}\subseteq
E$ such that $\{\hat{p}_{\delta }:\delta \in \hat{S}\}$ satisfies ($\dag $)
and ($\dag \dag $); in particular, $\hat{r}=r_{\delta }^{\hat{p}_{\delta
}(0)}$ for $\delta \in \hat{S}$ and $\eta _{\delta }^{\hat{p}_{\delta
}(0)}(n)$, $g^{\hat{p}_{\delta }(0)}(\delta ,n)=g(n)$ and $u^{\hat{p}%
_{\delta }(0)}(\delta ,n)$ are independent of $\delta $ for each $n<\hat{r}$%
. Moreover we can assume that there is $k\in \mathbb{Z}$ such that $k_{\delta
}=k$ for all $\delta \in \hat{S}$. Let $\hat{p}^{*}$ be the heart of $\{\hat{%
p}_{\delta }:\delta \in \hat{S}\}$ (so $\hat{p}^{*}\geq p^{*}\cup \{(i,%
\tilde{q})\}$).

Choose $m$ such that $2^{m}$ does not divide $c_{1}-c_{2}$. Let $M\geq \max
(\{g(n):n<\hat{r}\}\cup \{m\})$ and let 
\[
N=2^{1+(\hat{r}+1)M(d+1)}
\]
(where $d$ is the size of the domain of $\hat{p}^{*}-\{0\}$). Choose 
\[
\alpha _{0}<...<\alpha _{N-1}<\gamma <\delta _{0}<...<\delta _{N-1}
\]
where $\alpha _{j}\in S^{\prime }$, $\delta _{j}\in \hat{S}$, every ordinal
which occurs in $\hat{p}^{*}$ is $<\alpha _{0}$, and for all $j\leq N-1$
every ordinal which occurs in $p_{\alpha _{j}}$ or in $q_{\ell }^{\alpha
_{j}}$ ($\ell =1,2$) is less than $\alpha _{j+1}$ (where $\alpha _{N}$ is
taken to be $\gamma $); and for all $j<N-1$, every ordinal which occurs in $%
p_{\delta _{j}}$ is less than $\delta _{j+1}$. Then there is a condition $%
q_{0}\in Q_{0}$ which extends $\hat{p}^{*}(0)$ and each $p_{\alpha _{j}}(0)$
and $p_{\delta _{j}}(0)$ such that $q_{0}$ forces for all $j<N$: 
\[
\eta _{\delta _{j}}(\hat{r})=\gamma \text{; }\eta _{\delta _{j}}(\hat{r}%
+1)=\delta _{j-1}+1\text{; }u(\delta _{j},\hat{r})=y^{\alpha _{j}}\text{;
and }g(\delta _{j},\hat{r})=2^{m}
\]
where $\delta _{-1}=\gamma +1$.

As in the proof of Proposition \ref{zmt}, there is a condition $q^{\prime
}\in P_{i}$ and a subset $W^{\prime }$ of $N$ of size $\geq 2^{1+(\hat{r}%
+1)M}$ such that $q^{\prime }(0)=q_{0}$, $q^{\prime }\geq p_{\alpha _{j}}$
for all $j\leq N-1$, and $q^{\prime }\geq p_{\delta _{j}}\upharpoonright i$
for $j\in W^{\prime }$. Repeating the argument one more time and using the
facts that $q_{1}^{\alpha _{j}}$ and $q_{2}^{\alpha _{j}}$ force the same
value to $\varphi (y^{\alpha _{j}})$ and that $q_{1}^{\alpha
_{j}}\upharpoonright \alpha _{j}=q_{2}^{\alpha _{j}^{{}}}\upharpoonright
\alpha _{j}=\tilde{q}$, there is a subset $W=\{j,j_{o}\}$ of $W^{\prime }$
such that for any function $f:W\rightarrow \{1,2\}$ there is a condition $%
q_{*}\in P_{i+1}$ such that $q_{*}\upharpoonright i=q^{\prime }$ and $%
q_{*}(i)$ is an upper bound of $\{p_{\delta _{j}}(i),p_{\delta
_{j_{o}}}(i)\}\cup \{q_{1}^{\alpha _{j}},q_{2}^{\alpha _{j_{o}}}\}$. In a
generic extension $V[G^{\prime }]$ where $q_{*}\in G^{\prime }$ we have
(since $h(z_{\delta _{j},0}+K)=h(z_{\delta _{j_{o}},0}+K)=k$ and $g(\delta
,n)$ and $u(\delta ,n)$ are independent of $\delta \in \hat{S}$ for $n<\hat{r%
}$) that $2^{m}$ divides 
\[
h(u(\delta _{j},\hat{r}))-h(u(\delta _{j_{o}},\hat{r}))=h(y^{\alpha
_{j}}+K)-h(y^{\alpha _{j_{o}}}+K)=c_{1}-c_{2}
\]
which is a contradiction of the choice of $m$. This proves (III.3) and thus
finally completes the proof of Theorem \ref{free}.

\section{Theorem \ref{torsion}}

To prove Theorem \ref{torsion} we use a variant of the iterated forcing that
is described in section 1. Let $Q_{0}$ and $Q_{\psi }$ be as defined there.
We shall use a finite support iteration $P^{\prime }=\left\langle
P_{i}^{\prime },\dot{Q}_{i}:0\leq i<\omega _{2}\right\rangle $; the $\dot{Q}%
_{i}$ are defined inductively. We consider an enumeration, as before, of
names $\{\dot{\psi}_{i}:i<\omega _{2}\}$ for functions from $K$ to $\mathbb{Z} $%
. In $V^{P_{i}}$ we define

\[
\dot{Q}_{i}=\left\{ 
\begin{array}{ll}
\{0\} & \text{if the s.e.s }\mathcal{E}_{i}\text{ is torsion} \\ 
Q_{\dot{\psi}_{i}} & \text{otherwise}
\end{array}
\right. 
\]
We claim that if $G$ is $P^{\prime }$-generic, then in $V[G]$ (i) $\limfunc{%
Ext}(A,\mathbb{Z})$ is torsion and (ii) $\limfunc{Hom}(A,\mathbb{Z})=0$.

To see why (i) holds, consider $\psi \in \limfunc{Hom}(K,\mathbb{Z})$. For some 
$i\in \omega _{2}$, $\dot{\psi}_{i}$ is a name for $\psi $. In $V[G_{i}]$
either $\psi $ represents a torsion element of 
$\limfunc{Ext}(A,\mathbb{Z})$ 
or else, by construction, in $V[G_{i+1}]$ $\psi =\varphi |K$ for some $%
\varphi \in \limfunc{Hom}(F,\mathbb{Z})$, in which case $\psi $ represents the
zero element of $\limfunc{Ext}(A,\mathbb{Z})$.

To prove (ii), it suffices to show for all $i\in \omega _{2}$ that if $h\in 
\limfunc{Hom}(A,\mathbb{Z})^{V[G_{i+1}]}$, then $h\in \limfunc{Hom}(A,\mathbb{Z}%
)^{V[G_{i}]}$. If not, then $\dot{Q}_{i}\neq \{0\}$; but then by  
the arguments in section 4 it follows that 
$\mathcal{E}_{i}$ is torsion, so $\dot{Q}_{i}=\{0\}$,
 a contradiction.

\section{co-Moore spaces}

\smallskip \ Following \cite{HHS} we call a topological space $X$ a\textit{\
co-Moore space of type }$(G,n)$, where $n\geq 1$, if its reduced integral
cohomology groups satisfy 
\[
\tilde{H}^{i}(X,\mathbb{{Z)=}\left\{ 
\begin{array}{ll}
G & \text{if }i=n \\ 
0 & \text{otherwise.}
\end{array}
\right. } 
\]
For $n\geq 2$,  application of the Universal Coefficient Theorem
shows that

\begin{quote}
($\clubsuit $) there exist $B_{1}$ and $B_{2}$ such that 
$G\cong \limfunc{Hom%
}(B_{1},\mathbb{Z})\oplus \limfunc{Ext}(B_{2},\mathbb{Z)}$ where $\limfunc{Ext}%
(B_{1},\mathbb{Z})=0=\limfunc{Hom}(B_{2},\mathbb{Z)}$.
\end{quote}

Conversely, if $G$ satisfies ($\clubsuit $), then there is a co-Moore space
of type $(G,n)$ for any $n\geq 2$ (cf. \cite[Thm. 5]{HHS}, \cite{GG}%
). A sufficient condition for $G$ to be of the form ($\clubsuit $) is that $%
G=D\oplus C$ where $C $ is compact and
   $D  $ is isomorphic to a direct product
of copies of $\mathbb{Z}$ (\cite[Thm. 5]{HHS}). In a model of ZFC where every
W-group is free, this condition is necessary (cf. \cite[Thm. 3(a)]{HHS} and 
\cite[Thm. 2.20]{EH}); in particular the (torsion-free) rank of $C$ is of
the form $2^{\mu }$ for some infinite cardinal $\mu $. However, as a
consequence of our proofs we have:

\begin{corollary}
It is consistent with ZFC + $2^{\aleph _{0}}=2^{\aleph _{1}}=\aleph _{2}$
that there is a group $A$ of cardinality $\aleph _{1}$ such that 
$\limfunc{Hom}(A, \mathbb{Z}) =0$
but $\limfunc{Ext}(A,\mathbb{Z)}$ does not admit a compact topology.
\end{corollary}

\begin{corollary}
It is consistent with ZFC + $2^{\aleph _{0}}=2^{\aleph _{1}}=\aleph _{2}$
that for any $n\geq 2$ there is a co-Moore space of type $(F,n)$ where $F$
is the free abelian group of rank $\aleph _{2}$.
\end{corollary}

\begin{corollary}
\label{6.3}
It is consistent with ZFC + $2^{\aleph _{0}}=2^{\aleph _{1}}=\aleph _{2}$
that for any $n\geq 2$ there is a co-Moore space of type $(C,n)$ for some
uncountable torsion divisible group $C$.
\end{corollary}
 Compare Corollary \ref{6.3} with \cite[2.5 and 2.6]{GG}. The 
 conclusions of the corollaries are not provable in ZFC. Moreover, by an easy
modification we can replace $\aleph _{2}$ in the corollaries by any regular
cardinal greater than $\aleph _{1}$. (Note that by \cite[Thm. 5.6]{C1}, $%
\limfunc{Hom}(B_{1},\mathbb{Z})$ cannot be the free group of rank $\aleph _{1}$
if $\limfunc{Ext}(B_{1},\mathbb{Z})=0$.)

\end{document}